\documentclass[oneside,english,reqno]{amsart}
\usepackage{mathptmx}

\usepackage[T1]{fontenc}
\usepackage[latin9]{inputenc}
\usepackage{verbatim}
\usepackage{mathrsfs}
\usepackage{amsthm}
\usepackage{amssymb}
\usepackage{esint}
\usepackage{color}

\makeatletter
\numberwithin{equation}{section}
\numberwithin{figure}{section}
\theoremstyle{plain}
\newtheorem{thm}{\protect\theoremname}[section]
  \theoremstyle{plain}
  \newtheorem{lem}[thm]{\protect\lemmaname}
  \theoremstyle{definition}
  \newtheorem{defn}[thm]{\protect\definitionname}
    \theoremstyle{definition}
  \newtheorem{hyp}{Hypothesis}[section]
  \theoremstyle{remark}
  \newtheorem{rem}[thm]{\protect\remarkname}
  \theoremstyle{remark}
  \newtheorem{claim}[thm]{\protect\claimname}
  \theoremstyle{plain}
  \newtheorem{prop}[thm]{\protect\propositionname}

\makeatletter
\def\renewtheorem#1{%
  \expandafter\let\csname#1\endcsname\relax
  \expandafter\let\csname c@#1\endcsname\relax
  \gdef\renewtheorem@envname{#1}
  \renewtheorem@secpar
}
\def\renewtheorem@secpar{\@ifnextchar[{\renewtheorem@numberedlike}{\renewtheorem@nonumberedlike}}
\def\renewtheorem@numberedlike[#1]#2{\newtheorem{\renewtheorem@envname}[#1]{#2}}
\def\renewtheorem@nonumberedlike#1{  
\def\renewtheorem@caption{#1}
\edef\renewtheorem@nowithin{\noexpand\newtheorem{\renewtheorem@envname}{\renewtheorem@caption}}
\renewtheorem@thirdpar
}
\def\renewtheorem@thirdpar{\@ifnextchar[{\renewtheorem@within}{\renewtheorem@nowithin}}
\def\renewtheorem@within[#1]{\renewtheorem@nowithin[#1]}
\makeatother

  \providecommand{\claimname}{Hypothesis}
\DeclareMathOperator{\Lip}{Lip}

\makeatother

\usepackage{babel}
  \providecommand{\claimname}{Claim}
  \providecommand{\definitionname}{Definition}
  \providecommand{\lemmaname}{Lemma}
  \providecommand{\propositionname}{Proposition}
  \providecommand{\remarkname}{Remark}
\providecommand{\theoremname}{Theorem}

\begin{document}
\global\long\def\alg{\mathscr{F}}
\global\long\def\R{\mathbb{R}}
\global\long\def\N{\mathbb{N}}
\global\long\def\Z{\mathbb{Z}}
\global\long\def\C{\mathbb{C}}
\global\long\def\Q{\mathbb{Q}}
\global\long\def\F{\mathcal{F}}
\global\long\def\L{\mathscr{L}}
\global\long\def\E{\mathbb{E}}
\global\long\def\R{\mathbb{R}}
\global\long\def\H{\mathcal{H}}
\global\long\def\Pr{\mathbb{P}}
\global\long\def\borel{\mathfrak{B}}
\global\long\def\Norm{\mathscr{N}}
\global\long\def\Rd{\mathbb{R}^{d}}
\global\long\def\Hsd{\mathscr{H}}
\global\long\def\I{I_{\delta}}
\global\long\def\u{\mathbbm{1}}
\global\long\def\norm#1{\left\Vert #1\right\Vert }
\global\long\def\scal#1#2{\left\langle #1,#2\right\rangle _{H}}
\global\long\def\bh{\bar{h}}

\thanks{D.A. and M.M. are members of G.N.A.M.P.A. of the Italian Istituto Nazionale di Alta Matematica (INdAM)}

\title[Gradient contractivity in Wiener spaces]{Gradient contractivity of a rescaled resolvent on domains in Wiener spaces}

\author{Davide Addona, Giorgio Menegatti, Michele Miranda Jr}
\begin{abstract}
Given an abstract Wiener space $(X,\gamma,H)$, we consider an open set $O\subseteq X$ which satisfies certain smoothness and mean-curvature conditions. We prove that the rescaled resolvent operator associated to the Ornstein-Uhlenbeck operator with homogeneous Dirichlet boundary conditions on $O$ is gradient contractive in $L^p(X,\gamma)$ for every $p\in(1,\infty)$. This is the Gaussian counterpart of an analogous result for the rescaled resolvent operator associated to the Laplace operator $\Delta$ in $L^p$ with respect to the Lebesgue measure, $p\in[1,\infty)$, with homogeneous Dirichlet boundary conditions on a bounded convex open set $O\subseteq \R^n$.

\end{abstract}

\maketitle

\section{Introduction}

In this paper we consider an abstract Wiener space $(X,\gamma,H)$, where $X$ is a separable Banach space with a centered nondegenerate Gaussian measure $\gamma$ and $H$  its Cameron-Martin space $H$, and an open subset $O\subseteq X$ which satisfies suitable conditions. The aim of this paper is to prove that the rescaled resolvent $(J_\sigma^O=({\rm Id}-\sigma L_O)^{-1})_{\sigma>0}$ of the Ornstein-Uhlenbeck operator $L_O$ on $O$ with homogeneous Dirichlet boundary conditions satisfies
\begin{align}
\label{intr_res_contr_gaus}
\int_O\|\nabla_H J^O_\sigma f\|_{H}^pd\gamma\leq\int_O\|\nabla_Hf\|_H^pd\gamma, \quad f\in W^{1,q}_0(O,\gamma), \ \sigma>0, 
\end{align}
where $q=p$ if $p>1$ and $q>1$ if $p=1$. In particular, this implies that $(J_\sigma^O)_{\sigma>0}$ is $L^p$ gradient contractive for every $p\in(1,\infty)$.

This result extends to the Gaussian setting the analogous one proved in \cite[Appendix 1]{Bar}, where $O\subseteq \R^d$ is an open convex set with smooth boundary $\partial {O}$. In the quoted paper, the authors show that the rescaled resolvent associated to the Laplace operator on $O$ with homogeneous Dirichlet boundary conditions is $L^p$ gradient contractive  for every $p\in[1,\infty)$, i.e.,
\begin{align}
\label{intr_res_contr_leb}
\int_{O}\|\nabla J_\sigma f(x)\|_{\R^d}^pdx\leq \int_{O}\|\nabla f(x)\|_{\R^d}^pdx, \quad f\in W^{1,p}_0(O),
\end{align}
where $J_\sigma:=({\rm Id}-\sigma\Delta)^{-1}$, $\sigma>0$, is the rescaled resolvent of $\Delta$. We remark that for $p=\infty$ inequality \eqref{intr_res_contr_leb} has been proved in \cite{BrSt}. 

Let us compare formulae \eqref{intr_res_contr_gaus} and \eqref{intr_res_contr_leb}. In the latter, the operator $\nabla$ denotes the weak derivative in $L^p(O)$, while, in the former, the gradient $\nabla_H$ is the gradient along the directions of $H$. This is a typical issue of the infinite dimension: indeed, dealing with Gaussian measures in infinite dimension, an integration-by-parts formula, which is a crucial tool to prove the closability of the gradient operator, is verified if and only if one considers the gradient along $H$, which is usually denoted by $\nabla_H$. 

Another difference is the operator considered,  whose rescaled resolvent appears in the above inequalities: it is the Ornstein-Uhlenbeck operator in \eqref{intr_res_contr_gaus}, and the Laplace operator in \eqref{intr_res_contr_leb}. This follows from the fact that the Ornstein-Uhlenbeck operator plays the role of the Laplace operator when the Lebesgue measure is replaced by the Gaussian one. In finite dimension, this can be easily seen by taking into account an integration-by-parts formula: indeed, for smooth enough functions $f$ and $g$, we have
\begin{align*}
\int_{\R^d}\langle \nabla f,\nabla g\rangle_{\R^d}dx
= -\int_{\R^d}\Delta fg dx,
\end{align*}
and
\begin{align*}
\int_{\R^d}\langle \nabla f,\nabla g\rangle_{\R^d}d\gamma^d
= &\int_{\R^d}\langle \nabla f,\nabla g\rangle_{\R^d}e^{-|x|^2/2}dx
= -\int_{\R^d}L fg d\gamma^d,
\end{align*}
where $\gamma^d$ denotes the standard Gaussian measure on $\R^d$ and $Lf(x)=\Delta f(x)-\langle x,\nabla f(x)\rangle_{\R^d}$.

We further notice that, in the Gaussian context, gradient resolvent contractivity holds for every $p\in(1,\infty)$ (when $p=1$ we get the inequality not for functions $f\in W^{1,1}_0(O,\gamma)$ but for functions in smaller spaces), while \eqref{intr_res_contr_leb} is satisfied for every $p\in[1,\infty]$. The case $p=\infty$ is quite delicate in infinite dimension. Indeed, $W^{1,\infty}(X,\gamma)$ represents the intersection of $W^{1,p}(X,\gamma)$ with $p\in[1,\infty)$ and not the functions which have bounded $H$-gradient, and so we do not expect to extend \eqref{intr_res_contr_gaus} for $p=\infty$. On the contrary, for $p=1$ the main obstacle is that, to the best of our knowledge, a satisfactory theory of traces for Sobolev functions in infinite dimension is not available.
Traces of Sobolev functions at the boundaries of very smooth sets were considered for instance in \cite{add2,Cel}, while in \cite{Dap1, Dap} the Sobolev spaces of functions which "vanish" at the boundary are 
introduced to study maximal $L^2$ regularity for Dirichlet problems in infinite dimension. We also mention \cite{Hin4,Hin}, where the author deals with Sobolev spaces on domains.

The first attempt to provide a systematic study of traces of Sobolev functions on domains of abstract Wiener spaces appears in \cite{Cel}, where the authors give sufficient conditions on the domain $O$ to define a bounded operator ${\rm Tr}:W^{1,p}(O,\gamma)\rightarrow L^q(\partial O,\rho)$ for $p\in(1,\infty)$ and $q\in[1,p)$, where $\rho$ is the Hausdorff-Gauss surface measure of Feyel and de La Pradelle (see \cite{FeyDeL}). Under additional assumptions on $O$, it is proved that also $q=p$ can be achieved. Nothing can be said for the case $p=1$, and this shows that the gap between finite and infinite dimension is considerably big. We also stress that the case $p=1$ is not achieved either in \cite{add2} (which is strongly inspired by \cite{Cel}), whose characterization of Sobolev spaces $W^{1,p}_0(O,\gamma)$ of functions with null trace on $\partial O$ is widely used in this paper. We refer to \cite{Fer} for a theory of traces on domains in abstract Wiener spaces in $L^p$-spaces with respect to a weighted Gaussian measure, and to \cite{Add1,Bon} for an integration-by-parts formula on domains in Wiener spaces, which should be the starting point for a development in the study of traces by means of different techniques. Finally, for the case $p=1$, a possible alternative approach is to consider BV functions defined on open domains, which are investigated in \cite{Add1,bog2,LunMir}. However, it is still not clear how to extend to infinite dimension the theory of traces for BV functions in finite dimension.

To conclude, we spend few words on the assumptions on $O$ (see Hypotheses \ref{claim:regularity} and \ref{claim:Gaussian-curvature}). In \eqref{intr_res_contr_leb} the domain $O$ is assumed to be convex and with smooth boundary. Smoothness of the boundary and geometric properties for domains in infinite dimension are not easy to be defined, hence we translate the hypotheses on the cylindrical approximations $(O_n)_{n\in\N}$ of $O$, which can be considered as finite dimensional domains (see Hypothesis \ref{claim:regularity} and \cite{Dap}). As far as the regularity condition is concerned, we simply require for $\partial O_n$ the same smoothness of the domain $O$ in \cite{Bar}. Further, the convexity of $O$ implies that the mean curvature at every point $x\in \partial  O$ is non-negative, and this fact is crucial to prove \eqref{intr_res_contr_leb}. Since $\partial O$ is smooth, it follows that the mean curvature at $x\in\partial O$ is the divergence of the outer normal to $O$ at $x$. In the Gaussian setting, the classical divergence operator div is replaced by the Gaussian divergence ${\rm div}_\gamma$, which on smooth vector fields $F$ acts as ${\rm div}_\gamma F(x)={\rm div}F(x)-\langle x,F(x)\rangle$. Hence, it seems to be reasonable to ask that for every $x\in \partial O_n$ the Gaussian divergence of the normal to $O_n$ at $x$ is non-negative. This is indeed the right choice, and for every $x\in \partial O_n$, this value is defined as the {\it Gaussian mean curvature} at $x\in\partial O_n$. So, under the assumptions that $\partial O_n$ is smooth and that the Gaussian mean curvature is non-negative at $x$ for every $x\in\partial O_n$, definitely with respect to $n$, we are able to prove \eqref{intr_res_contr_gaus}. 

The paper is organized as follows. In Section \ref{sec:setting} we provide the classical basic results on abstract Wiener spaces. Further, we state the assumptions on the domain $O$, define the Sobolev spaces $W^{1,p}(O,\gamma)$, $W^{1,p}_0(O,\gamma)$ and the Ornstein-Uhlenbeck operator $L_O$ with homogeneous Dirichlet boundary conditions on $\partial O$ by means of the theory of Dirichlet forms (see for instance \cite{davies,Fuku80}), and recall the main results of \cite{add2} which will be used in the paper. In Section \ref{sec:res_contr_fin_dim1} we show \eqref{intr_res_contr_gaus} when $X=\R^d$ for some $d\in \N$. In Section \ref{sec:res_contr_infin_dim} we extend the results of Section \ref{sec:res_contr_fin_dim1} when $X$ is a separable Banach space. To this aim, we split this section into two parts. In the former we prove \eqref{intr_res_contr_gaus} when $O$ is a cylindrical domain with respect to a fixed orthonormal basis $\{h_n:n\in\N\}$ of $H$ (see Hypothesis \ref{hyp:dom_cyl_case}) and its finite dimensional projection $\mathcal O$ is a domain with smooth boundary and with non-negative Gaussian mean curvature at every point of its boundary. In the latter we show that \eqref{intr_res_contr_gaus} also holds true for non-cylindrical domains $O$ under suitable assumptions on the cylindrical approximations $(O_n)_{n\in\N}$ of $O$. Finally, in Section \ref{sec:examples} we provide some examples of domains $O$ in abstract Wiener spaces which satisfy our conditions.

\subsection{Notation}
Given a separable Banach space $X$ and its topological dual $X^*$, we denote by $\|\cdot\|_X$ its norm and by $\langle \cdot,\cdot\rangle_{X\times X^*}$ its duality. 

Let $A$ be an open set in $\R^{d}$. For every $k\in\N\cup\{\infty\}$ we denote by $C^{k}(A)$ the set of functions on $A$ which
are $k$-times differentiable on $A$ with continuous derivatives up to order $k$. $C_b^{k}(\R^{d})$
is the set of $k$ times differentiable functions which are continuous and bounded together with their derivatives up to order $k$. $C_{c}^{\infty}(\R^{d})$ is the subspace of $C_{b}^{\infty}(\R^{d})$ of functions with compact support.

Let $\alpha\in(0,1)$. We denote by $C^{\alpha}(A)$ the set of $\alpha$-H\"older continuous functions on $A$. We denote by $C_{\rm loc}^{\alpha}(A)$ the set of functions $f$ which are $\alpha$-H\"older continuous on every bounded open set $U\subseteq A$.
We denote by $C^{2,\alpha}(A)$ the set of functions $f\in C^{2}(A)$ such that $D^mf\in C^\alpha(A)$ for every multi-index $m$ with length $|m|=2$.
We denote by $C_{\rm loc}^{2,\alpha}(A)$ the set of functions $f\in C^{2}(A)$ such that $D^mf\in C_{\rm loc}^\alpha(A)$ for every multi-index $m$ with length $|m|=2$.

If $A$ is replaced by $\overline A$ in the definition of the above spaces, we mean that the functions have a continuous extension up to $\overline A$.

For every $O\subseteq \R^d$ with non-empty interior and every $k\in \N\cup\{\infty\}$, we denote by $C_{0}^k(O)$ the subset of $C^k(\R^d)$ of functions which vanish out of an open set $A\subseteq O$ with positive distance from $O^c$, and we denote by $C_{0}^k(\overline O)$ the subset of $C^k(\R^d)$ consisting of functions which vanish on $O^c$. We denote by $C_c^k(O)$ the subset of $C^k(O)$ of functions which have compact support in $O$.

If a function $f$ is defined on $A\subseteq X$, then $\overline f$ denotes the trivial extension of $f$ to the whole space $X$, i.e., $\overline f=f$ on $A$ and $\overline f=0$ on $A^c$.

A subset of $\R^d$ is said to be $C^{2,\alpha}$-regular if its boundary is locally a graph of a function in $C^{2,\alpha}(\R^{d-1})$.


Let $H$ be a separable Hilbert space. We denote by $\langle\cdot,\cdot,\rangle_H$ its inner product and by $\mathcal L(H)$ the space of linear bounded operator $A:H\rightarrow H$ endowed with the norm-operator topology. 

We say that a nonnegative operator $A\in \mathcal L(H)$ is nuclear (or trace-class) in $H$ if there exists an orthonormal basis $\{e_n:n\in\N\}$ of $H$ such that
\begin{align*}
{\rm Tr} (A):=\sum_{n\in\N}\langle A e_n,e_n\rangle_H<\infty.
\end{align*}
We stress that ${\rm Tr} (A)$ does not depend on choice the basis $\{e_n:n\in\N\}$. We denote by $\mathcal L^+_1(H)$ the space of trace-class operators on $H$.

We say that the operator $A\in\mathcal L(H)$ is a Hilbert-Schmidt operator if there exists an orthonormal basis $\{e_n:n\in\N\}$ of $H$ such that
\begin{align*}
\|A\|_{\mathscr L_2(H)}^2:=\sum_{n=1}^\infty\|Ae_n\|_H^2<\infty.
\end{align*}
The above series does not depend on the choice of the basis, and we denote by $\mathscr L_2(H)$ the subspace of $\mathcal L(H)$ consisting of Hilbert-Schmidt operators. The space $(\mathscr L_2(H),\|\cdot\|_{\mathscr L_2(H)})$ is a separable Hilbert space. If $A\in\mathscr L_2(H)$ then $ AA^*\in \mathcal L^+_1(H)$.

Let $Y$ be a separable Banach space; we denote by $L^{p}(X,\gamma,Y)$ as the space of (the equivalence classes of) Bochner integrable functions $F:X\to Y$ such that 
\begin{align}
\label{norma_Lp}
\|F\|_{L^{p}(X,\gamma,Y)}=\bigg(\int_{X}\|F\|_{Y}^{p}\ d\gamma\bigg)^{1/p}<\infty,
\end{align}
see e.g. \cite{Die}. This space, endowed with the norm \eqref{norma_Lp}, is a Banach space. 

With $\borel(X)$ we denote the Borel $\sigma$-algebra of $X$, and $B_X(x,r)$, $x\in X$, $r>0$, denotes the open ball of $X$ with center $x$ and radius $r$.

\section{Preliminary results}
\label{sec:setting}

\subsection{\label{sub:Fundamentals-about-Wiener}Fundamentals about abstract Wiener spaces}

Let us recall some definitions and properties of Gaussian measures on separable Banach spaces. For a detailed treatment we refer to the monograph  \cite{Bog}.

Let $X$ be a separable Banach space and let $\gamma$ be a centered nondegenerate Gaussian measure on $X$. Since $\gamma$ is a Gaussian measure, the elements of $X^*$ can be seen as elements of $L^2(X,\gamma)$. We consider the embedding $j:X^{*}\hookrightarrow L^{2}(X,\gamma)$ and the reproducing kernel $X_{\gamma}^{*}$ is defined as the closure of $j(X^{*})$ in $L^{2}(X,\gamma)$. $(X^*_\gamma,\|\cdot\|_{L^2(X,\gamma)})$ is  a separable Hilbert space, and we define $R_{\gamma}:X_{\gamma}^{*}\rightarrow(X^{*})^{*}$ as 
\[
R_{\gamma}(f)(g)=\int_{X}fg\ d\gamma,
\qquad f\in X_{\gamma}^{*}, \ \in X^{*}.
\]
$R_{\gamma}$ has range in $X$ and it is injective. We define the Cameron-Martin space $H$ as $R_{\gamma}(X_{\gamma}^{*})\subseteq X$, and it inherits a structure of separable Hilbert space from $X^*_\gamma$ through $R_\gamma$: for every $h\in H$ we denote by $\hat h\in X^*_\gamma$ the unique element such that $R_\gamma \hat h=h$ and $\langle h, k\rangle _H:=\langle \hat h,\hat k\rangle_{L^2(X,\gamma)}$ for every $h,k\in H$. The space $H$ is continuously and densely embedded into $X$, and we denote by $c_H$ the smaller positive constant $c$ which satisfies 
\begin{align}
\label{cost_C_H}
\|h\|_X\leq c\|h\|_H, \qquad h\in H.
\end{align}
We introduce the operator $Q:X^*\rightarrow H$, defined as $Q(x^*)=R_\gamma (j(x^*))\in X$ for every $x^*\in X^*$, and we fix an orthonormal basis $\{h_{n}:n\in\N\}$ of $H$ consisting of elements of $Q(X^*)$, i.e.,  $h_{n}= Q(x^*_n)$, with $x^*_n\in X^*$, for every $n\in\N$. Let us fix $n\in\N$ and set
\begin{align*}
F_{n}:={\rm Span}\{ h_{1},\ldots,h_{n}\}, \quad \pi_{n}(x)=\sum_{i=1}^{n}\langle x,x^*_i\rangle_{X\times X^*} h_{i}, \quad x\in X.
\end{align*}
Hence, $\pi_{n}$ is the projection of $X$ on the finite dimensional subspace $F_{n}$ of $H$. For every $n\in\N$ we denote by $\gamma_{n}$ the image measure of $\gamma$ by means of $\pi_{n}$, i.e., $\gamma_{n}:=\gamma\circ\pi_{n}^{-1}$. The measure $\gamma_{n}$ is a nondegenerate centered Gaussian measure on $F_n$. If we set $X_{n}^{\perp}:={\rm Ker}(\pi_{n})$, then it follows that $X=F_n\oplus X_{n}^{\perp}$, and for every $n\in\N$ we introduce the isomorphism $\Pi_{n}:F_n\rightarrow \R^n$ defined as 
\begin{align*}
\Pi_{n}x:=(x_1,\ldots,x_n), 
\qquad x=\sum_{i=1}^n x_ih_i.
\end{align*}

For every $k\in\N\cup\{\infty\}$ we introduce the set $\mathcal{F}C_{b}^k(X)$ of bounded cylindrical functions which are $k$-times Fr\'echet differentiable, i.e., the functions $f:X\to \R$ such that there exist $n\in\N$, $l^*_{1},\ldots,l^*_{n}\in X^*$ and $\varphi\in C_b^{k}(\R^{n})$ such that $f(x)=\varphi(\langle x,l^*_{1}\rangle_{X\times X^*}, \ldots, \langle x,l^*_{n}\rangle_{X\times X^*})$ for every $x\in X$. We stress that for every $k\in\N\cup\{\infty\}$, $\mathcal{F}C_{b}^{k}(X)$ is dense in $L^{p}(X,\gamma)$ for every $p\in[1,\infty)$.




\subsection{$H$-derivative and Sobolev spaces}

For every $h\in H$ and $f\in\mathcal{F}C_{b}^{\infty}(X)$ we define the $H$-derivative $\partial_{h}f:X\to \R$ of $f$ along $h$ as 
\[
\partial_{h}f(x)=\lim_{\varepsilon\rightarrow0}\frac{f(x+\varepsilon h)-f(x)}{\varepsilon}, \quad x\in X,
\]
and its formal adjoint
$\partial_{h}^{*}f=\partial_{h}f-f\widehat{h}$.

For $f\in\mathcal{F}C _{b}^{\infty}(X)$, there exists a unique $\nabla_{H}f:X\rightarrow H$, called $H$-gradient of $f$, such that
\[
\partial_{h}f(x)=\left\langle \nabla_{H}f(x),h\right\rangle {}_{H}, \quad x\in X, \ h\in H.
\]
If $f$ is a Lipschitz function on $X$, then $\partial_hf$ can be defined $\gamma$-a.e. and the essentially bounded function $\nabla_{H}f$ can be defined $\gamma$-a.e. and identified with an element of $L^{\infty}(X,\gamma,H)$ (see e.g. \cite[Theorem 5.11.2]{Bog}). 

We say that $f:X\rightarrow \R$ is $H$-differentiable at $x\in X$ if there exists $F(x)\in H$ such that
\begin{align*}
f(x+h)=f(x)+\langle F(x),h\rangle_H+o(\|h\|_H), \quad \|h\|_H\rightarrow 0.
\end{align*}
We set $\nabla_Hf(x):=F(x)$ for every $x\in X$ such that $f$ is $H$-differentiable at $x$.

For $f:X\rightarrow H$, we say that $f$ is $H$-differentiable at $x\in X$ if there exists a Hilbert-Schmidt operator $D_{H}f(x)$ on $H$ such that
\begin{align*}
f(x+h)=f(x)+D_Hf(x)h+o(\|h\|_H), \quad \|h\|_H\rightarrow0.
\end{align*}
$D_{H}f(x)$ is said $H$-derivative of $f$ at $x$, and for every $h\in H$ we have
\[
D_{H}f(x)(h) =\lim_{\varepsilon\rightarrow0}\frac{f(x+\varepsilon h)-f(x)}{\varepsilon},
\]

%
%
%

Let $f:X\rightarrow\R$ be such that $\nabla_{H}f$ is defined on the whole $X$. We say that $f$ is twice $H$-differentiable at $x\in X$ if $f$ is $H$-differentiable and $\nabla_{H}f(x)$ has $H$-derivative, which we denote by $D_{H}^{2}f(x)$. 

The set of smooth cylindrical vector-valued functions $\mathcal{F}C_{b}^{\infty}(X,H)$, defined as the linear span of the functions $\phi h$ where $\phi\in\mathcal{F}C_{b}^{\infty}(X)$
and $h\in H$. Let $f\in\mathcal{F}C_{b}^{\infty}(X,H)$, $x\in X$ and $h\in H$, then $D_Hf(x)$ and $\partial_{h}f(x)$ are well-defined and $\partial_{h}f(x)=D_{H}f(x)(h)$.

The integration-by-parts formula
\begin{align*}
\int_X \partial_hfd\gamma=\int_X f\hat hd\gamma, \quad f\in \mathcal{F}C _b^\infty(X), \ h\in H,
\end{align*}
is the key tool to prove that for every $p\in[1,\infty)$, the operators $\nabla_H:\mathcal{F}C _b^\infty(X)\rightarrow L^p(X,\gamma,H)$ and $(\nabla_H,D^2_H):\mathcal{F}C _b^\infty(X)\rightarrow L^p(X,\gamma,H)\times L^p(X,\gamma,\mathscr L_2(H))$ are closable in $L^p(X,\gamma)$, and the operator $D_H:\mathcal{F}C _b^\infty(X)\rightarrow L^p(X,\gamma,\mathscr L_2(H))$  is closable in $L^p(X,\gamma,H)$. We still denote by $\nabla_H$, $D_H$ and $D^2_H$, respectively, the closure of these operators in $L^p$.

\begin{defn}
We denote by $W^{1,p}(X,\gamma)$ the domain of $\nabla_H$ in $L^p (X,\gamma)$, by $W^{1,p}(X,\gamma,H)$ the domain of $D_H$ in $L^p (X,\gamma,H)$ and by $W^{2,p}(X,\gamma)$ the domain of $(\nabla_H, D^2_H)$ in $L^p (X,\gamma)$. These spaces are Banach spaces if endowed with the norms
\begin{align*}
& \|f\|_{W^{1,p}(X,\gamma)}  :=\left(\|f\|_{L^p(X,\gamma)}^p+\|\nabla_Hf\|_{L^p(X,\gamma,H)}^p\right)^{1/p}, \quad f\in W^{1,p}(X,\gamma), \\
& \|f\|_{W^{2,p}(X,\gamma)}  :=\left(\|f\|_{L^p(X,\gamma)}^p+\|\nabla_Hf\|_{L^p(X,\gamma,H)}^p+\|D^2_Hf\|_{L^p(X,\gamma,\mathcal H_2(H))}^p\right)^{1/p}\!\!\!\!\!\!, \quad f\in W^{2,p}(X,\gamma), \\
& \|f\|_{W^{1,p}(X,\gamma,H)}  :=\left(\|f\|_{L^p(X,\gamma,H)}^p+\|D_Hf\|_{L^p(X,\gamma,\mathcal H_2(H))}^p\right)^{1/p}, \quad f\in W^{1,p}(X,\gamma,H),
\end{align*}
respectively. Finally, if $p=2$ the above spaces are Hilbert spaces.
\end{defn}

Since $\nabla_H$ is closed and densely defined $L^2(X,\gamma)$, its adjoint operator ${\rm div}_\gamma:=\nabla_H^*$ is closed, densely defined and satisfies
\begin{eqnarray}
\int_{X}\left\langle f,\nabla_{H}g\right\rangle _{H}\ d\gamma =  -\int_{X}\mbox{div}_{\gamma}fg\ d\gamma,
\label{eq:divergence-1}
\end{eqnarray}
for every $g\in W^{1,2}(X,\gamma)$ and every $f$ in the domain $D({\rm div}_\gamma)\subseteq L^2(X,\gamma,H)$ of ${\rm div}_\gamma$. From \cite[Theorem 5.8.2]{Bog} it follows that $W^{1,2}(X,\gamma,H)\subseteq D({\rm div}_\gamma)$ and
\begin{align*}
{\rm div}_\gamma F=\sum_{n=1}^\infty\left( \partial_{h_n} F_n- F_n\widehat h_n\right), \quad F\in W^{1,2}(X,\gamma,H),
\end{align*}
where the series converges in $L^2(X,\gamma)$ and $F_n:=\langle F,h_n\rangle_H$. Further, the above formula is independent of the choice of the basis of $H$ and $\|{\rm div}_\gamma F\|_{L^2(X,\gamma)}\leq \|F\|_{W^{1,2}(X,\gamma,H)}$.
In particular, if $F\in \mathcal{F}C _b^\infty(X,H)$ we get
\begin{equation}
\mbox{div}_{\gamma}F=\sum_{i=1}^{m}\partial_{k_{i}}F_{i}-\sum_{i=1}^{m}\widehat{k_{i}}F_{i},
\label{eq:divergence}
\end{equation}
where $F$ satisfies $F=\sum_{i=1}^m F_ik_i$, with $F_i\in \mathcal{F}C _b^\infty(X)$ for every $i=1,\ldots,m$ and $\{k_1,\ldots,k_m\}$ are orthonormal in $H$.

\subsection{The Sobolev spaces $W^{1,p}(O,\gamma)$ and $W_{0}^{1,p}(O,\gamma)$}
\label{sub:sobolev_spaces}


Let $O\subseteq X$ be an open set. We denote by $\Lip(O)$ the set of Lipschitz functions on $O$, and by ${\rm Lip}_c(O)$ the subset of ${\rm Lip}(O)$ whose elements vanish out of an open set $A$ with positive distance from $O^{c}$. For every $m\in\N$ we denote by ${\rm Lip}_{c,m}(O,H)$ the set of $H$-valued Lipschitz functions on $O$ of the form
\begin{align*}
f=\sum_{i=1}^mf_ih_i, \qquad f_i\in {\rm Lip}_{c}(O), \ h_i\in H, \ i=1,\ldots,m.
\end{align*}
The set
\begin{align*}
{\mathcal F}{\rm Lip}_c(O,H):=\bigcup_{m\in\N}{\rm Lip}_{c,m}(O;H),
\end{align*}
is dense in $L^p(O,\gamma,H)$. Finally, for every $m\in\N$ and $f\in{\rm Lip}_{c,m}(O,H)$, the function ${\rm div}_\gamma \overline f$ is defined $\gamma$-almost everywhere as
\begin{align*}
{\rm div}_\gamma  \overline f=\sum_{i=1}^m\left(\partial_{h_i} \overline f_i- \overline f_i\widehat h_i\right)
\end{align*}
and ${\rm div}_\gamma \overline f\in L^p(X,\gamma)$ for every $p\in(1,\infty)$.

The proof of the following lemma can be found in \cite[Lemma 2.1]{Add}. 
\begin{lem}
\label{lem:Dirichlet_SObolev_O}For every $p\in[1,\infty)$, the
operator $\nabla_{H}:\Lip(O)\rightarrow L^{p}(O,\gamma, H)$ is closable in $L^p(O,\gamma)$. We still denote by $\nabla_H$ its closure.

\end{lem}

\begin{defn}
\label{def:Sobolev }(Sobolev Spaces) Let $p\in[1,\infty)$.
$W^{1,p}(O,\gamma)$ is the domain of $\nabla_{H}$ in $L^p(O,\gamma)$. The space $W^{1,p}(O,\gamma)$ is a Banach space if endowed with the norm
\[
\left\Vert f\right\Vert _{W^{1,p}(O,\gamma)}:=\left(\left\Vert f\right\Vert _{L^{p}(O,\gamma)}^p+\left\Vert \nabla_{H}f\right\Vert _{L^{p}(O,\gamma,H)}^p\right)^{1/p}, \quad f\in W^{1,p}(O,\gamma),
\]
and $W^{1,2}(O,\gamma)$ is a Hilbert space with inner product 
\[
\left\langle f,g\right\rangle _{W^{1,2}(O,\gamma)}=\left\langle f,g\right\rangle _{L^{2}(O,\gamma)}+\left\langle \nabla_{H}f,\nabla_{H}g\right\rangle _{L^{2}(O,\gamma,H)},\quad  f,g\in W^{1,2}(O,\gamma).
\]
\end{defn}

\begin{defn}
\label{def:C_b_H}
We denote by $\mathcal H^{1}(X)$ the set of all continuous functions $f$ (not necessarily bounded) which are $H$-differentiable on $X$ and such that $\nabla_{H}f$
is bounded and continuous with values in $H$.

$\mathcal H_{b,0}^{1}(O)$ is the subset of $\mathcal H^1(X)$ of bounded functions $f$ which vanishes out of an open set $A$ with positive distance from $O^{c}$. 
\end{defn}
\begin{rem}
\label{rem:spazio_funz_test_prop}
$\mathcal H_{b,0}^{1}(O)$ is not empty. Taking advantage from the results in \cite{Su98} In \cite[Lemma 2.2]{add2} it has been proved that the subset of $\mathcal H^1(X)$ whose elements vanish out of an open set $A$ with positive distance from $O^c$ is not empty. We simply remark that the function $F_{B,\varepsilon}$ provided in the quoted lemma is also bounded, and so it belongs to $\mathcal H^1_{b,0}(X)$.
\end{rem}

The bounded elements of $\mathcal H^1(X)$, and so in particular the elements of $\mathcal H_{b,0}^1(X)$, can be approximated in a useful way, as the next result shows.

\begin{lem}
\label{lem:convergence_projec_1} If $f\in \mathcal H^1(X)$ is bounded, then $f\circ\pi_{n}\rightarrow f $ in ${W^{1,p}(X,\gamma)}$ as $n\rightarrow\infty$ for every $p\in[1,\infty)$.
\end{lem}
\begin{proof}
We set $f_{n}=f\circ\pi_{n}$ for every $n\in\N$. From \cite[Corollary 3.5.8]{Bog} it follows that $\pi_{n}x\rightarrow x$ as $n\rightarrow\infty$ in $X$ for $\gamma$-a.e. $x\in X$. The continuity of $f$ implies that $f_{n}\rightarrow f$
$\gamma$-a.e. in $X$, and so $f_{n}$ converges to $f$ in $L^{p}(X,\gamma)$ by the dominated convergence theorem. Moreover,
\[
\nabla_{H}f_{n}(x)=\pi_{n}(\nabla_{H}f(\pi_n(x))), \qquad x\in X, \ n\in\N,
\]
and by the definition of $\pi_{n}$ it follows that $\|\pi_{n}(h)\|_{H}\leq\|h\|_{H}$ and $\pi_{n}(h)\rightarrow h$ in $H$ as $n\rightarrow\infty$ for every $h\in H$. Hence,
\begin{align*}
\|\nabla_{H}f_{n}(x)-\nabla_{H}f(x)\|_{H}
\leq & \|\nabla_{H}f_n(x)-\pi_{n}(\nabla_{H}f(x))\|_{H}
+\|\pi_{n}\circ\nabla_{H}f(x)-\nabla_{H}f(x)\|_{H} \\
\leq & \|\nabla_{H}f(\pi_{n}(x))-\nabla_{H}f(x)\|_{H}+\|\pi_{n}(\nabla_{H}f(x))-\nabla_{H}f(x)\|_{H}.
\end{align*}
The two addends in the very last right hand-side of the above chain of inequalities vanish as $n\rightarrow\infty$, the former for $\gamma$-a.e. $x\in X$ since $\nabla_Hf$ is continuous, and the latter for every $x\in X$ due to the convergence of $\pi_{n}$ in $H$. Recalling that $\nabla _Hf$ is bounded, from the dominated convergence theorem, we get the thesis.
\end{proof}


\begin{defn}
\label{def:Sobolev_Dir} 
For every $p\in[1,\infty)$, we denote by $W_{0}^{1,p}(O,\gamma)$ the closure of ${\rm Lip}_c(O)$ in $W^{1,p}(O,\gamma)$. 
\end{defn}

The proof of the following lemma can be obtained repeating verbatim that of \cite[Lemma 2.3]{add2}, hence we omit it.
\begin{lem}
\label{lem:dens_H_funct_sob_dir}
The closure of $\mathcal H_{b,0}^1(O)$ in $W^{1,p}(O,\gamma)$ coincides with $W^{1,p}_0(O,\gamma)$ for every $p\in[1,\infty)$.
\end{lem}

\subsection{The Ornstein-Uhlenbeck semigroup and operator}
For every $f\in C_b(X)$ and every $t>0$ we set
\begin{align*}
T_tf(x):=\int_X f(e^{-t}x+\sqrt{1-e^{-2t}}y)d\gamma(y), \quad x\in X,
\end{align*}
and $T_0f=f$. Thanks to the equality
\begin{align*}
\int_XT_tfd\gamma=\int_Xfd\gamma, \qquad f\in C_b(X), \ t\geq 0, 
\end{align*}
for every $t\geq 0$ and every $p\in[1,\infty)$ the operator $T_t$ extends to a contraction operator on $L^p(X,\gamma)$ still denoted by $T_t$. It turns out that $(T_t)_{t\geq0}$ is a strongly continuous semigroup of contractions on $L^p(X,\gamma)$ and we denote by $L_p$ its infinitesimal generator. If $p=2$ we write $L$ instead of $L_2$. Further, $\gamma$ is the unique invariant measure for $(T_t)_{t\geq0}$, $\mathcal{F}C _b^\infty(X)$ is a core for the domain $D(L_p)$ of $L_p$, $D(L_p)=W^{2,p}(X,\gamma)$ for every $p\in(1,\infty)$ and
\begin{align*}
\int_XLfgd\gamma=-\int_X\langle \nabla_Hf, \nabla_Hg\rangle_Hd\gamma, 
\end{align*}
for every $f\in D(L)$ and $g\in W^{1,2}(X,\gamma)$ (see \cite[Chapter 5]{Bog}). In particular, if $f\in\mathcal{F}C _b^2(X)$ then $f\in D(L_p)$ for every $p\in(1,\infty)$ and
\begin{align*}
Lf=\sum_{i=1}^n\partial_{k_ik_i}^2f-\sum_{i=1}^n \hat k_i \partial_{k_i}f,
\end{align*}
where $f(x)=\varphi(\hat k_1(x),\ldots, \hat k_n(x))$ for every $x\in X$, with $\varphi\in C^2(\R^n)$ and $k_1,\ldots,k_n\in Q(X^*)$ are orthonormal vectors in $H$.

Now we define the Ornstein-Uhlenbeck operator on $O$, starting from the bilinear form
\begin{align*}
a(f,g):=\int_O\langle \nabla_Hf,\nabla_Hg\rangle_Hd\gamma, \qquad f,g\in W^{1,2}_0(O,\gamma).    
\end{align*}
From the theory of Dirichlet forms (see \cite{davies}), there exists a unique closed operator $L_O$ with dense domain $D(L_O)\subseteq W_{0}^{1,2}(O,\gamma)$ in $L^2(O,\gamma)$ such that
\begin{align}
\label{def_OU_operator_forme}
-\int_{O}\left\langle \nabla_{H}f,\nabla_{H}g\right\rangle _{H}\ d\gamma=\int_{O}L_Of\cdot g\ d\gamma, \qquad f\in D(L_O), \ g\in W^{1,2}_0(O,\gamma).
\end{align}

\begin{defn}
\label{def:Ornstein} 
$L_O$ is called Ornstein-Uhlenbeck operator with homogeneous Dirichlet boundary conditions on $\partial O$.
We set 
\begin{align*}
J^O_\sigma:=({\rm Id}-\sigma L_O)^{-1}=\sigma^{-1}(\sigma^{-1}{\rm Id}-L_O)^{-1}, \qquad \sigma>0.
\end{align*}
$J_\sigma$ is a bounded operator on $L^{2}(X,\gamma)$ with range equals to $D(L_O)$ and the family $(J^O_\sigma)_{\sigma>0}$ is called {\textit{rescaled resolvent} of $L_O$}. For every $\sigma>0$ we denote by $G_\sigma^O$ the resolvent of $L_O$, i.e., $G_\sigma^O=(\sigma I-L_O)^{-1}=\sigma^{-1} J_{\sigma^{-1}}^O$ for every $\sigma>0$.

 If $O=X$ we simply write $L$ instead of $L_X$.
\end{defn}

\begin{rem}
From the theory of symmetric Markov semigroups (see \cite[Section 1.4]{davies}), the semigroup $(T^O_2(t))_{t\geq0}$ associated to $L_O$ in $L^2(O,\gamma)$  extends from $L^\infty(O,\gamma)$ to a positive contraction strongly continuous semigroup $(T^O_p(t))_{t\geq0}$ on $L^p(O,\gamma)$ for every $p\in[1,\infty)$. These semigroups are consistent in the sense that, if $1\leq p\leq q<\infty$, if $f\in L^q(O,\gamma)$ then $T^O_p(t)f=T^O_q(t)f$ for every $t\geq0$. For every $p\in[1,\infty)$ we denote the infinitesimal generator of $(T^O_p(t))_{t\geq0}$ by $L_{O,p}$, and if $p=2$ we simply write $L_O$ instead of $L_{O,2}$. For every $\sigma>0$ we denote by $G_\sigma^{O,p}$ and by $J^{O,p}_\sigma$ the resolvent and the rescaled resolvent of $L_{O,p}$, respectively, and we recall that both $G_\sigma^{O,p}$ and $J^{O,p}_\sigma$ are continuous linear operators on $L^p(O,\gamma)$.
\end{rem}

\subsection{The case when $O$ is the sublevel set of a function $G$}
Hereafter we assume the following (see \cite[Hypothesis 3.1]{add2} and \cite[Hypothesis 3.1]{Cel}).
\begin{claim}
\label{claim:regularity}
Let $G:X\rightarrow \R$ and $\delta>0$ be such that:
\begin{enumerate} 
\item[i)]$G\in \mathscr H^{1}(X)$;
\item[ii)] $\nabla_{H}G$ is everywhere $H$-differentiable (in particular, $G$ is twice $H$-differentiable), with
derivative $D_{H}^{2}G$ and $\|D_{H}^{2}G\|_{\mathscr L_2(H)}$ uniformly bounded;
\item[iii)] $G^{-1}(0)\neq\emptyset$;
\item[iv)] $\|\nabla_{H}G\|_{H}^{-1}\in L^{\infty}(X)$;
\item[v)] $LG$ is bounded on $G^{-1}(-\delta,\delta)$.
\end{enumerate}
Under the above assumptions, we set that $O:=G^{-1}((-\infty,0))$ and, for every $n\in\N$, $G_n:=G\circ \pi_n$ and $O_n:=(G_n)^{-1}((-\infty,0))$. In particular, $x\in O_n$ if and only if $\pi_n(x)\in O$. 
\end{claim}
Without loss of generality we may assume that $\overline O\neq X$. Indeed, if $\overline O=X$, then $W^{1,p}_0(O,\gamma)=W^{1,p}(X,\gamma)$ since $\gamma(\partial O)=0$ from \cite[Remark 3.2]{add2}.

\begin{rem}
\label{rmk:hyp_G_sodd}
If $G$ fulfills Hypothesis \ref{claim:regularity}, then it satisfies \cite[Hypothesis 3.1]{add2}.    
\end{rem}

\begin{lem}
\label{lem:H_O_n_da_H_O}
Let $f\in \mathscr H^1_{b,0}(O)$. Then, the function $f_n:=f\circ \pi_n$ belongs to $\mathcal H_{b,0}^1(O_n)$ for every $n\in\N$.   
\end{lem}
\begin{proof}
It is enough to prove that $f_n$ vanishes on $A_n^c$, where $A_n$ is an open set with positive distance from $O_n^c$. To prove this fact, let $A\subseteq O$ be an open set with positive distance $d$ from $O^c$ such that $f$ vanishes on $A^c$, and let us fix $n\in\N$. We define $A_n:=\{x\in X:\pi_n(x)\in A\}$. 
We claim that $A_n$ is open, has positive distance from $O_n^c$ and that $f_n$ vanishes on $A_n^c$. The last assertion is easy to prove. Indeed, for every $x\in A_n^c$ it follows that $\pi_n x\in A^c$, and so $f_n(x)=f(\pi_n x)=0$. It remains to show that $A_n$ is open and has positive distance from $O_n^c$. Let $x\in A_n$ and let $z\in O_n^c$. Then, $\pi_n(x)\in A$, $\pi_n (z)\in O^c$ and
\begin{align}
\label{stima_dist_proj}
d\leq \|\pi_n(x)-\pi_n(z)\|_X
\leq & \sum_{i=1}^n|\langle x-z,x^*_i\rangle_{X\times X^*}\|h_i\|_X
\leq \overline c\|x-z\|_X \sum_{i=1}^n\|x^*_i\|_{X^*},
\end{align}
where $\overline c$ has been introduced in Subsection \ref{sub:Fundamentals-about-Wiener}. This implies that $A_n$ has positive distance from $O_n^c$. Finally, arguing as in \eqref{stima_dist_proj}, for every $x,y\in X$ we get
\begin{align*}
\|\pi_n(y)-\pi_n(x)\|_X\leq \overline c \|x-y\|_X\sum_{i=1}^n\|x^*_i\|_{X^*}.   
\end{align*}
Let $x\in A_n$. Hence, $\pi_n(x)\in A$ and if $\delta>0$ fulfills $B_X(\pi_n(X),\delta)\subseteq A$, then from the above computations we get
\begin{align*}
B_{X}\left(x,\left(\overline c\sum_{i=1}^n\|x^*_i\|_{X^*}\right)^{-1}\delta\right)\subseteq A.
\end{align*}
We conclude that $A_n$ is open, which gives $f_n\in \mathcal H_{b,0}^1(O_n)$. 
\end{proof}


From Remark \ref{rmk:hyp_G_sodd} and \cite[Theorem 4.1]{add2} we state the following result.
\begin{thm}
\label{thm:approximation:W_0}
The following are equivalent:
\begin{enumerate}
\item[i)]$f\in W_{0}^{1,p}(O,\gamma)$;
\item[ii)]the trivial extension $\overline f$ of $f$ out of $O$ belongs to $W^{1,p}(X,\gamma)$. 
\end{enumerate}
\end{thm}

\begin{rem}
\label{rem:grad_est_nulla}
For every $f\in W_0^{1,p}(O,\gamma)$, we have $\overline {\nabla_Hf}=\nabla_H\overline f$. The fact is obvious if $f\in{\rm Lip}_c(O)$ (or $f\in \mathcal H_{b,0}^1(O)$), the general case follows by approximation.
\end{rem}

\begin{rem}
\label{rem:lower_semicont}
Let $p\in(1,\infty)$ and $r$ with $1< r\leq p$, we have that the function on $W^{1,p}_0(O,\gamma)$ given by $y\mapsto \|\nabla_Hy\|_{L^p(O,\gamma,H)}$ is lower semi-continuous with respect to the topology of $L^r(O,\gamma)$. To prove this fact, we recall that if $q$ is the conjugate exponent of $p$ then $L^q(O,\gamma,H)$ is the dual of $L^p(O,\gamma,H)$, see for instance \cite[Chapter 4, Theorem 1]{Die}. Therefore, for every $f\in W^{1,p}_0(O,\gamma)$, from Theorem \ref{thm:approximation:W_0} and Remark \ref{rem:grad_est_nulla} we get
\begin{align*}
\|\nabla_Hf\|_{L^p(O,\gamma,H)}
= & \sup_{G\in L^q(O,\gamma,H), \ \|G\|_{L^q(O,\gamma,H)}\leq 1}\int_O\langle \nabla_Hf,G\rangle_Hd\gamma \\
= & \sup_{G\in \mathcal F{\rm Lip}_{c}(O,H), \ \|G\|_{L^q(O,\gamma,H)}\leq 1}\int_O\langle \nabla_H f, G\rangle_Hd\gamma \\
= & \sup_{G\in \mathcal F{\rm Lip}_{c}(O,H), \ \|G\|_{L^q(O,\gamma,H)}\leq 1}\int_X\langle \nabla_H\overline f,\overline G\rangle_Hd\gamma \\
= & \sup_{G\in \mathcal F{\rm Lip}_{c}(O,H), \ \|G\|_{L^q(O,\gamma,H)}\leq 1}\int_X\overline f\ {\rm div}_\gamma\overline Gd\gamma \\
= & \sup_{G\in \mathcal F{\rm Lip}_{c}(O,H), \ \|G\|_{L^q(O,\gamma,H)}\leq 1}\int_O f\ {\rm div}_\gamma  Gd\gamma.
\end{align*}
Since for every $G\in \mathcal F{\rm Lip}_{c}(O,H)$ the map 
\begin{align*}
f\mapsto \int_Of{\rm div}_\gamma  Gd\gamma    
\end{align*}
is continuous with respect to the topology of $L^r(O,\gamma)$, the lower semicontinuity of $y\mapsto \|\nabla_Hy\|_{L^p(O,\gamma,H)}$ with respect to the topology of $L^r(O,\gamma)$ follows from the lower semicontinuity of the supremum.
\end{rem}

We conclude this section by showing that, if $f\in \mathscr H^{1}_{b,0}(O)$ and we set $f_n=f\circ \pi_n$ for every $n\in\N$, then $(\overline{J_\sigma^{O_n}f_n})_{n\in\N}$ converges to $\overline{J_{\sigma}^Of}$ in $W^{1,2}(X,\gamma)$ as $n$ goes to infinity.
\begin{prop}
\label{prop:conv_ris_n}
Let $f\in \mathscr H^{1}_{b,0}(O)$ and let us set $f_n=f\circ \pi_n$ for every $n\in\N$. Then, 
$(\overline{J_\sigma^{O_n}(f_n)})_{n\in\N}$ converges to $\overline{J_{\sigma}^O(f)}$ in $W^{1,2}(X,\gamma)$ as $n$ goes to infinity.    
\end{prop}
\begin{proof}
The fact that $f_n\in\mathscr H^1_{b,0}(O_n)$ for every $n\in\N$ follows from Lemma \ref{lem:H_O_n_da_H_O}. Further, from Remark \ref{rem:grad_est_nulla} and Theorem \ref{thm:approximation:W_0} it follows that $\overline {J^O_\sigma(f)}\in W^{1,2}(X,\gamma)$
 and $\overline {\nabla_HJ^O_\sigma(f)}=\nabla_H\overline {J^O_\sigma(f)}$ for every $\sigma>0$.
 
We fix $\sigma>0$ and $n\in\N$, and we consider the trivial extension $\overline{J_\sigma^{O_n}(f_n)}$ of $J_\sigma^{O_n}(f_n)\in W^{1,2}_0(O_n,\gamma)$. From Remark \ref{rem:grad_est_nulla} and Theorem \ref{thm:approximation:W_0} we infer that $\overline{J_\sigma^{O_n}(f_n)}\in W^{1,2}(X,\gamma)$ and $\nabla_H \overline{J_\sigma^{O_n}(f_n)}=\overline {\nabla_HJ_\sigma^{O_n}(f_n)}\in W^{1,2}(X,\gamma)$. Further, since 
\begin{align}
\label{sol_deb_L_n}
J_\sigma^{O_n}(f_n)-\sigma L_{O_n}J_\sigma^{O_n}(f_n)=f_n
\end{align}
on $O_n$, multiplying both the sides of this equation by $J_\sigma^{O_n}(f_n)$ and integrating on $O_n$ with respect to $\gamma$ we get
\begin{align*}
\|J_\sigma^{O_n}(f_n)\|_{L^2(O_n,\gamma)}^2
+\sigma \|\nabla_H J_\sigma^{O_n}(f_n)\|_{L^2(O_n,\gamma,H)}^2 = & \int_{O_n} J_\sigma^{O_n}(f_n) f_n d\gamma \\
\leq & \|J_\sigma^{O_n}(f_n)\|_{L^2(O_n,\gamma)}\|f_n\|_{L^2(O_n,\gamma)},
\end{align*}
which gives 
\begin{align*}
\|\overline{J_\sigma^{O_n}(f_n)}\|_{L^2(X,\gamma)}^2
+\sigma \|\nabla_H \overline{J_\sigma^{O_n}(f_n)}\|_{L^2(X,\gamma,H)}^2 \leq & \|\overline{J_\sigma^{O_n}(f_n)}\|_{L^2(X,\gamma)}\|\overline{f_n}\|_{L^2(X,\gamma)},    
\end{align*}
from which it follows that
$\|\overline{J_\sigma^{O_n}(f_n)}\|_{L^2(X,\gamma)}\leq \|\overline{f_n}\|_{L^2(X,\gamma)}$. This implies that 
\begin{align}
\label{stima_norma_sob_res_n}
\|\overline{J_\sigma^{O_n}(f_n)}\|_{L^2(X,\gamma)}^2
+\sigma \|\nabla_H \overline{J_\sigma^{O_n}(f_n)}\|_{L^2(X,\gamma,H)}^2 \leq & \|\overline{f_n}\|_{L^2(X,\gamma)}^2
\end{align}
and, since $(\overline {f_n})$ converges in $L^2(X,\gamma)$ (see Lemma \ref{lem:convergence_projec_1}), we infer that $(\overline{J_\sigma^{O_n}(f_n)})_{n\in\N}$ is a bounded sequence in $W^{1,2}(X,\gamma)$. Therefore, up to a subsequence, $(\overline{J_\sigma^{O_n}(f_n)})_{n\in\N}$ weakly converges to some $u$ in $W^{1,2}(X,\gamma)$ and
\begin{align*}
\|u\|_{W^{1,2}(X,\gamma)}^2
\leq \liminf_{n\to\infty}\|\overline{J_\sigma^{O_n}(f_n)}\|_{W^{1,2}(X,\gamma)}^2.
\end{align*}

We claim that there exists $v\in W^{1,2}(X,\gamma)$ with $v=0$ $\gamma$-a.e. in $O^c$ and $u=v$ $\gamma$-a.e. in $X$. If the claim is true, then from Theorem \ref{thm:approximation:W_0} it follows that $u_{|O}=v_{|O}\in W^{1,2}_0(O,\gamma)$.

\noindent Let $g\in{\rm Lip}_c(\overline{O}^c)$, and let us consider the sequence $(g_n)_{n\in\N}$ defined as $g_n=g\circ\pi_n$ for every $n\in\N$. Arguing as in Lemma \ref{lem:H_O_n_da_H_O}, it is possible to prove that $g_n\in{\rm Lip}_c((\overline {O}^c)_n)$ for every $n\in\N$, where $(\overline O^c)_n:=\{x\in X:\pi_n(x)\in \overline O^c\}$. In particular, if $B\subseteq \overline O^c$ is an open set which has positive distance from $\overline O$ and $g=0$ on $B^c$,  then $B_n=\{x\in X:\pi_n(x)\in B\}$ is an open set with positive distance from $(\overline{O}^c)_n$ and $g_n=0$ on $B_n^c$. Hence, we get
\begin{align*}
\int_Xugd\gamma=\lim_{n\to\infty}\int_X\overline {J_\sigma^{O_n}(f_n)}g_nd\gamma=0,    
\end{align*}
since the supports of $\overline {J_\sigma^{O_n}(f_n)}$ (which is the set $\overline{O_n}$) and of $g_n$ (which is the set $\overline{B_n}$) have positive distance. Indeed, if $d:={\rm dist}(\overline B,\overline O)$, then for every $x\in B_n$ we get $\|\pi_n(x)-z\|_X\geq d$ for every $z\in \overline O$. In particular, for every $y\in O_n$ we get $\pi_n(y)\in O$, and so $\|\pi_n(x)-\pi_n(y)\|_X\geq d$ for every $x\in B_n$ and every $y\in O_n$. Arguing as in \eqref{stima_dist_proj} we infer that there exists a positive constant $\widetilde d$ such that $\widetilde d\leq \|x-y\|_X$ for every $x\in B_n$ and every $y\in O_n$, which implies that $\overline O_n$ and $\overline B_n$ have positive distance.

The arbitrariness of $g$ and the density of ${\rm Lip}_c(\overline{O}^c)$ in $L^2(\overline{O}^c,\gamma)$ give $u=0$ $\gamma$-a.e. in $\overline O^c$. Finally, since under our assumptions we have $\gamma(\partial O)=0$ (see \cite[Remark 3.2]{add2}), we conclude that the function $v$, defined as $v=u$ on $X\setminus \partial O$ and $v=0$ on $\partial O$, fulfills $v=0$ $\gamma$-a.e. in $O^c$ and $v=u$ $\gamma$-a.e. in $X$. The claim is so proved.

Let us show that $u_{|O}=J_\sigma^{O}(f)$. For every $g\in \mathscr H^1_{b,0}(O)$,we set $g_n:=g\circ \pi_n$ for every $n\in\N$. From Lemma \ref{lem:H_O_n_da_H_O} we infer that $g_n\in \mathscr H^1_{b,0}(O_n)$ for every $n\in\N$ and multiplying both the sides of \eqref{sol_deb_L_n} by $g_n$ we get
\begin{align}
\label{int_parti_u_3}
\sigma\int_{O_n}\langle \nabla_H{J_\sigma^{O_n}(f_n)},\nabla_H g_n\rangle_Hd\gamma
=\int_{O_n}J_\sigma^{O_n}(f_n) g_nd\gamma-\int_{O_n}f_n g_nd\gamma.
\end{align}
From Lemma \ref{lem:convergence_projec_1} we know that $(g_n)_{n\in\N}$ converges to $\overline g$ in $W^{1,2}(X,\gamma)$ as $n$ tends to $\infty$. Hence, the weak convergence of $(\overline{J_\sigma^{O_n}(f_n)})_{n\in\N}$, the strong convergence of $(g_n)_{n\in\N}$ and Remark \ref{rem:grad_est_nulla} give
\begin{align}
\sigma \int_O\langle \nabla_H u,\nabla_Hg\rangle_Hd\gamma
= & \sigma\int_X\langle \nabla_Hu,\overline {\nabla_Hg}\rangle_{H}d\gamma 
=  \sigma\int_X\langle \nabla_Hu,\nabla_H \overline g\rangle_Hd\gamma \notag \\
= & \sigma\lim_{n\to\infty}\int_X\langle \nabla_H\overline{J_\sigma^{O_n}(f_n)},\nabla_H g_n\rangle_Hd\gamma \notag \\
= & \sigma\lim_{n\to\infty}\int_{O_n}\langle \nabla_H{J_\sigma^{O_n}(f_n)},\nabla_Hg_n\rangle_Hd\gamma \notag \\
= & \lim_{n\to\infty}
\left(\int_{X}\overline {J_\sigma^{O_n}(f_n)}g_nd\gamma-\int_{X}\overline{f_n}\overline g_nd\gamma\right) \notag  \\
= & \int_Xugd\gamma-\int_X\overline fgd\gamma \notag \\
= & \int_Ougd\gamma-\int_Ofgd\gamma. 
\label{int_parti_u_4}
\end{align}

From \eqref{int_parti_u_4} we get
\begin{align}
\sigma\int_O\langle \nabla_H u,\nabla_Hg\rangle_H d\gamma=\int_Ougd\gamma-\int_Ofgd\gamma,     
\label{int_parti_u_5}
\end{align}
and the arbitrariness of $g\in\mathscr H^1_{b,0}(O)$ implies that $u_{|O}=J_\sigma^O(f)$, i.e., $u=\overline{J_\sigma^O(f)}$ $\gamma$-a.e. 

It remains to prove that $(\overline{J_\sigma^{O_n}(f_n)})_{n\in \N}$ converges to $u$ in $W^{1,2}(X,\gamma)$. To this aim, we show that $(\overline{J_\sigma^{O_n}(f_n)})_{n\in \N}$ converges to $u$ in ${W}^{1,2}(X,\gamma)$ with respect to the equivalent norm
\begin{align*}
\|f\|_{\sigma}^2
= \|f\|_{L^2(X,\gamma)}^2+\sigma\|\nabla_Hf\|_{L^2(X,\gamma,H)}^2, \qquad f\in\mathcal FC_b^1(X).
\end{align*}
Multiplying both the sides of \eqref{sol_deb_L_n} by $J^{O_n}_\sigma(f_n)$, integrating on $X$ and letting $n$ go to infinity, since $(\overline{f_n})_{n\in\N}$ converges to $\overline f$ in $L^2(X,\gamma)$ as $n$ tends to infinity, we get
\begin{align}
\lim_{n\to\infty}\|\overline{J_\sigma^{O_n}(f_n)}\|_{\sigma}^2
= & 
\lim_{n\to\infty}\left(\|\overline{J_\sigma^{O_n}(f_n)}\|_{L^2(O_n,\gamma)}^2
+\sigma \|\nabla_H \overline{J_\sigma^{O_n}(f_n)}\|_{L^2(O_n,\gamma,H)}^2\right) \notag \\
= & \int_O f u d\gamma.
\label{lim_ris_n_spazio_sob_tilde}
\end{align}

Let us consider a sequence $(u_m)_{m\in\N}\subseteq \mathscr H^1_{b,0}(O)$ which converges to $u_{|O}$ in $W^{1,2}_0(O,\gamma)$ and, for every $n\in\N$, let us set $u_{m,n}=u_m\circ \pi_n$. Hence, $u_{m,n}\in \mathscr H^1_{b,0}(O_n)$ for every $n\in\N$. Multiplying both the sides of \eqref{sol_deb_L_n} by $u_{m,n}$, integrating on $X$ with respect to $\gamma$ and applying the definition of $L_{O_n}$, we infer that
\begin{align}
\int_X\overline{J_{\sigma}^{O_n}(f_n)}u_{m,n}d\gamma+\sigma\int_X\langle \nabla_H\overline {J_{\sigma}^
{O_n}(f_n)},\nabla_Hu_{m,n}\rangle_Hd\gamma= \int_X\overline f_n u_{m,n}d\gamma
\label{int_parti_u_finale_1}
\end{align}
for every $m,n\in\N$. From Lemma \ref{lem:convergence_projec_1} we know that the sequence $(u_{m,n})_{n\in\N}$ converges to $\overline u_m$ in $W^{1,2}(X,\gamma)$ as $n$ goes to infinity. Therefore, letting $n$ go to infinity in \eqref{int_parti_u_finale_1} we infer that
\begin{align*}
\int_Xu\overline u_{m}d\gamma+\sigma\int_X\langle \nabla_Hu,\nabla_H\overline u_m\rangle_Hd\gamma
=\int_X\overline f\overline u_md\gamma, \qquad m\in\N.    \end{align*}
Recalling that $u_m=0$ on $O^c$, we infer that
\begin{align}
\int_Ouu_{m}d\gamma+\sigma\int_O\langle \nabla_Hu,\nabla_Hu_m\rangle_Hd\gamma
=\int_O fu_md\gamma        
\label{int_parti_u_finale_2}
\end{align}
for every $m\in\N$. Letting $m$ go to infinity in \eqref{int_parti_u_finale_2} and recalling that $(u_m)_{m\in\N}$ converges to $u_{|_O}$ in $W^{1,2}_0(O,\gamma)$ as $m$ goes to infinity, it follows that
\begin{align}
\|u\|^2_{\sigma}
=\|u\|_{L^2(O,\gamma)}^2+\sigma\|\nabla_Hu\|_{L^2(O,\gamma,H)}^2=\int_O f ud\gamma.    
\label{norma_u_spazio_sob_tilde}
\end{align}
From \eqref{lim_ris_n_spazio_sob_tilde} and \eqref{norma_u_spazio_sob_tilde} we infer that
\begin{align*}
\|u\|^2_{\sigma}
= \lim_{n\to\infty}\|\overline{J_\sigma^{O_n}(f_n)}\|_{\sigma}^2,
\end{align*}
which combined with the weak convergence of $(\overline{J_\sigma^{O_n}(f_n)})_{n\in \N}$ to $u$ in $W^{1,2}(X,\gamma)$ gives the strong convergence of $(\overline{J_\sigma^{O_n}(f_n)})_{n\in \N}$ to $u$ in $W^{1,2}(X,\gamma)$.

Finally, we have shown that every subsequence of $(\overline{J_\sigma^{O_n}(f_n)})_{n\in \N}$ admits a subsequence which converges to $\overline{J_\sigma^O(f)}$ in $W^{1,2}(X,\gamma)$ as $n$ goes to infinity, which is enough to conclude that the whole sequence $(\overline{J_\sigma^{O_n}(f_n)})_{n\in \N}$ converges to $u$ in $W^{1,2}(X,\gamma)$.
\end{proof}

\section{Resolvent contractivity in finite dimension}
\label{sec:res_contr_fin_dim1}

In this section we consider the finite dimensional case $X=\R^{d}$ endowed with the standard Gaussian measure $\gamma=\gamma^d$, with density $\theta_d(x)=(2\pi)^{-d/2}e^{-|x|^2/2}$ for every $x\in\R^d$. In this case, $H=X$ and the $H$-inner product is the Euclidean inner product. 

Let $O\subseteq\R^{d}$ be an open set which is $C^{2,\alpha}$-regular. The aim of this section is to prove \eqref{intr_res_contr_gaus} in this finite dimensional setting. 

\begin{rem}
\label{rem:regularity1}
For every $u\in C^2_0(\overline O)$ we get $L_Ou(x)=\overline Lu(x):=\Delta u(x) -\langle x,\nabla u(x)\rangle_{\R^d}$ for every $x\in O$. Since $I-\sigma \overline {L}$ is an elliptic operator with regular coefficients for every $\sigma>0$, it follows that for every bounded $C^{2,\alpha}$-regular domain $\Omega$ and every $f\in C^{\alpha}(\overline{\Omega})$, there exists a unique solution $g\in C^{2,\alpha}(\overline \Omega)$ of $g-\sigma\overline Lg=f$ on $\overline \Omega$ with $g_{|\partial\Omega}=0$ (see e.g. \cite[Theorem 6.14]{Gil}).
\end{rem}

From Remark \ref{rem:regularity1} we infer that if $y\in C^\infty_c( O)$ then $u=J_\sigma(y)$ belongs to $C_0^{2}(\overline O)$ (in the sense of the representatives) and it satisfies $(I-\sigma\overline {L})u=y$. Further,  from elliptic regularity (see e.g. \cite[Section 6.3.1, Theorem 3]{Eva}) we infer that for any bounded open subset $U\subseteq O$ the function $u$ belongs to $C^\infty(U)$. In particular, $u\in C^\infty(O)$.

Let $y\in C_c^{\infty}(O)$ and let $u\in C^\infty(O)\cap C^2_0(\overline O)$ satisfy $u-\sigma \overline Lu=u-\sigma L_Ou=y$.
We introduce the functions $\varphi$ and $\varphi_{\varepsilon}$ as follows:
\begin{align}
\varphi(x)=\|\nabla u(x)\|_{\R^d},\ \ \ \varphi_{\varepsilon}(x)=\sqrt{\varepsilon^{2}+\|\nabla u(x)\|_{\R^d}^{2}}.
\label{def_varphi}
\end{align}
Let us notice that $\varphi_\varepsilon\in C^\infty(O)$ for every $\varepsilon>0$. We prove the following result, which is the equivalent of \cite[Lemma 8.2]{Bar}.
\begin{lem}
\label{lem: phi1}
For every $x\in O$, we have
\begin{equation}
\frac{(\varphi(x))^{2}}{\varphi_{\varepsilon}(x)}-\sigma \overline L\varphi_{\varepsilon}(x)<\|\nabla y(x)\|_{\R^d}.
\label{eq:lemma1}
\end{equation}
\end{lem}
\begin{proof}
For every $j\in\{1,\ldots,d\}$ we have 
\begin{equation*}
\varphi_{\varepsilon}D_{j}\varphi_{\varepsilon}=\sum_{i=1}^{d}D_{i}uD_{ij}^{2}u,
\end{equation*}
which gives
\begin{align}
\|\nabla\varphi_{\varepsilon}\|_{\R^d}^{2}
= & \frac{\displaystyle\sum_{j=1}^{d}\left(\displaystyle \sum_{i=1}^{d}D_{i}uD_{ij}^{2}u\right)^{2}}{\varphi_{\varepsilon}^{2}}\leq\frac{\displaystyle\sum_{j=1}^{d}\left(\displaystyle\sum_{i=1}^{d}(D_{i}u)^{2}\displaystyle\sum_{i=1}^{d}(D_{ij}^{2}u)^{2}\right)}{\varphi_{\varepsilon}^{2}} \notag \\
= & \frac{\displaystyle\sum_{i=1}^{d}\left(D_{i}u\right)^{2}}{\varphi_{\varepsilon}^{2}}\sum_{i,j=1}^{d}\left(D_{ij}^{2}u\right)^{2}=\frac{\varphi^{2}}{\varphi_{\varepsilon}^{2}}\sum_{i,j=1}^{d}\left(D_{ij}^{2}u\right){}^{2}<\sum_{i,j=1}^{d}\left(D_{ij}^{2}u\right){}^{2}.\label{eq:inequality1}
\end{align}
Since on smooth functions $f$ the operator $\overline L$ reads as $\overline Lf( x)=\Delta f( x)-\langle x,\nabla f( x)\rangle_{\R^d}$ for every $ x\in O$, for every $i=1,\ldots,d$ we get
\begin{align*}
D_i(\overline L f)( x)
= & D_i(\Delta f)( x)-D_i(\langle  x, \nabla f( x)\rangle_{\R^d}) = \Delta (D_i)f( x)-\langle  x, \nabla (D_if)( x)\rangle_{\R^d}-D_if( x) \\
= & \overline L(D_i)f( x)-D_if( x), \quad  x\in O,
\end{align*}
whence
\begin{align}
\label{derOU_OUder1}
D_i(\overline Lf)( x) D_if( x)
= & \overline L(D_if)( x)D_if( x)-(D_if( x))^2 
\leq  \overline L(D_if)( x)D_if( x)
\end{align}
for every $ x\in O$ and $i=1,\ldots,d$.
Indeed, from the definition of $\varphi_\varepsilon$ and $\overline L$ we get
\begin{align*}
\varphi_\varepsilon( x) \sigma \overline L\varphi_\varepsilon( x)
= & \varphi_\varepsilon( x) \sigma \sum_{i=1}^d \left(D_{ii}\varphi_\varepsilon( x)- x_iD_i \varphi_\varepsilon(x)\right)  \\
= & \sigma\left(-\|\nabla \varphi_\varepsilon(x)\|_{\R^d}^2+\!\!\sum_{i,j=1}^d|D^2_{ij} u(x)|^2+\!\!\sum_{i,j=1}^d(D^3_{iij} u(x)- x_iD^2_{ij}u( x))D_ju( x)\right) \\
= & \sigma\left(-\|\nabla \varphi_\varepsilon(x)\|_{\R^d}^2+\sum_{i,j=1}^d|D^2_{ij} u(x)|^2+\sum_{i,j=1}^d\overline L(D_ju)( x)D_ju( x)\right)
\end{align*}
for every $ x\in O$. We recall that $\|\nabla \varphi_\varepsilon\|_{\R^d}^2\leq \sum_{i,j=1}^d(D^2_{ij} u)^2$ (see \eqref{eq:inequality1}), hence
\begin{align*}
\varphi_\varepsilon( x) \sigma \overline L\varphi_\varepsilon( x)
\geq & \sigma \sum_{j=1}^d\overline L(D_ju)( x)D_ju( x)
\end{align*}
for every $ x\in O$, and from \eqref{derOU_OUder1} we infer that
\begin{align}
\varphi_\varepsilon( x) \sigma \overline L\varphi_\varepsilon( x)
\geq & \sigma  \sum_{j=1}^d D_j (\overline Lu)( x)D_ju( x)=   \sum_{j=1}^d D_j (\sigma \overline Lu)( x)D_ju( x) \notag \\
= & \sum_{j=1}^d (D_j u( x)-D_jy( x))D_ju( x) =\|\nabla u( x)\|_{\R^d}^2- \langle \nabla u( x) \nabla y( x)\rangle_{\R^d} \notag \\
\geq & \varphi( x)^2-\varphi( x)\|\nabla y( x)\|_{\R^d},
\label{dis_varphi_epsilon_finale1}
\end{align}
where we have used the fact that $\sigma \overline Lu=u-y$. The thesis follows by dividing the first and the last side in \eqref{dis_varphi_epsilon_finale1} by $\varphi_\varepsilon$ and rearranging the terms.
\end{proof}

\subsection{Gaussian curvature}
We use the spaces $W^{1,p}(O,\gamma)$
and $W_{0}^{1,p}(O,\gamma)$ introduced in Section \ref{sec:setting}. We stress that, since we are in finite dimension, $C_{c}^{\infty}(O)$ is dense in $W_{0}^{1,p}(O,\gamma)$. We define $L_O$ and, for every $\sigma>0$, $J_O^\sigma$, accordingly to Definition \ref{def_OU_operator_forme}.

Let $x\in\partial O$ and let $\nu$ be the outer normal to $\partial O$ at $x$. By the $C^{2,\alpha}$-regularity of $O$, there exist a neighbourhood $U$ of $x$, an open set $V\subseteq \R^{d-1}$ and a smooth function $\psi:V\to \R$ such that, up to rotation which gives $\nu=-\xi_d$, $U\cap \partial O=\{(\xi',\psi(\xi')):\xi'\in V\}$, i.e., $\partial O$ is locally a graph of a function $\psi$ with the vertical axis oriented inside $O$, and we define the (inner) mean curvature of $\partial O$ at $x\in\partial O$ as $\Delta\psi(\xi'_0)$, where $\xi'_0\in V$ satisfies $x=(\xi'_0,\psi(\xi'_0))$ since $\nabla \psi(\xi'_0)=0$. Equivalently, if $O$ is the sublevel of a $C^{2}$-function $g$ with $\nabla g\neq 0$ on $\partial O$, we have that the mean curvature at $x\in \partial O$ is 
\begin{equation}
H_{\partial O}(x)=\frac{\Delta g(x)}{\|\nabla g(x)\|_{\R^d}}-\frac{\langle D^2g(x) \nabla g(x),\nabla g(x)\rangle_{\R^d} }{\|\nabla g(x)\|_{\R^d}^{3}},\label{eq:mean curvature1}
\end{equation}
where $D^2g$ is the Hessian matrix of $g$. 
We point out that the geometric mean curvature is \eqref{eq:mean curvature1} multiplied by $(d-1)^{-1}$.

\begin{defn}
\label{def:gaussian_curvature1}
If $O\subseteq\R^{d}$ is $C^{2,\alpha}$-regular, the (inner) Gaussian curvature at $x\in\partial O$ is $H^\gamma_{\partial O}(x)=H_{\partial O}(x)-\langle x,\nu(x)\rangle_{\R^d}$, where $H_{\partial O}$ is the mean curvature and $\nu$ is the outer normal to $\partial O$. 
\end{defn}

\begin{hyp}
\label{hyp:gauss_curv1}
We assume that the open set $O$ is $C^{2,\alpha}$-regular for some $\alpha>0$ and that it has non negative (inner) Gaussian curvature.
\end{hyp}

As we said, given a point $\overline x\in\partial O$ there exists a neighbourhood $U$ of $\overline x$ such that $U\cap\partial O$ can be seen as the graph of a smooth function $\psi:V\subseteq \R^{d-1}\to \R$.  In $\R^d$, we consider the rotation $R^{-1}$ (change of coordinates) $x\mapsto \xi$ (i.e., $R\xi=x$)  centered at $\overline x$ such that $R^{-1}(U\cap \partial O)$ is actually the graph of $\psi$. In the new system of coordinates, the outer normal to $\partial O$ at $\overline x$ is oriented as $-\xi_d$, where in this new system of coordinates $\xi_d$ is the $d$-th element of the basis. In particular, $(\xi',\psi(\xi'))=U\cap \partial O$ for every $\xi'\in V$, $\overline \xi=(\overline \xi',\psi(\overline \xi'))=\overline x$ with $\overline \xi'\in V$ and $\nabla_{\xi'}\psi(\overline \xi')=0$.


We notice that the operator $\overline L$ introduced in Remark \ref{rem:regularity1} is invariant under the new system of coordinates $\xi$, i.e., $\overline Lw(\xi)=\overline Lu(x)$ where $u,w\in C^\infty(\R^d)$ and $w(\xi)=w(R^{-1}x)=u(x)$ for every $\xi\in\R^d$ and $x=R\xi$. Indeed, since $x=R\xi$ and $R$ is an orthogonal matrix, for $u(x)=u(R\xi):=w(\xi)=w(R^{-1}x)$ it follows that
\begin{align*}
D_{x_i}u(x)=\langle \nabla_\xi w(\xi),(R^{-1})_{ i}\rangle_{\R^d}, \quad D^2_{x_ix_i}u(x)= \langle D^2_\xi w(\xi)(R^{-1})_{ i}, (R^{-1})_{ i}\rangle_{\R^d}, 
\end{align*}
where $D_\xi w$ is the Hessian matrix of $w$ with respect to the variable $\xi$ and $(R^{-1})_i$, $i=1,\ldots,d$, denotes the $i$-th row of the matrix $R^{-1}$. Hence,
\begin{align*}
\overline Lu(x)
= & \Delta_x u(x)-\langle x,\nabla_x u(x)\rangle_{\R^d}
= \sum_{i,j,k=1}^d D^2_{\xi_j\xi_k}w(\xi) (R^{-1})_{ji} (R^{-1})_{ki}
-\sum_{i,j=1}^d x_i D_{\xi_j}w(\xi)(R^{-1})_{ji},
\end{align*}
and recalling that $\sum_{i=1}^d  (R^{-1})_{ji} (R^{-1})_{ki}=\delta_{jk}$ and that $\sum_{i=1}^d(R^{-1})_{ji}x_i=\xi_j$, it follows that
\begin{align*}
\overline Lu(x)
= \Delta_\xi w(\xi)-\langle \xi,\nabla_\xi w(\xi)\rangle_{\R^d}=\overline Lw(\xi).
\end{align*}
Further, if $\varphi$ and $\varphi_\varepsilon$ are the functions introduced in \eqref{def_varphi}, it follows that
\begin{align*}
\varphi(x)^2
= &\sum_{i=1}^d\sum_{j,k=1}^mD_{\xi_j}w(\xi)(R^{-1})_{j i}D_{\xi_k}w(\xi)(R^{-1})_{k i} \\
= &\sum_{j,k=1}^dD_{\xi_j}w(\xi)D_{\xi_k}w(\xi)\delta_{jk} 
= \|\nabla_\xi w(\xi)\|_{\R^d}^2,
\end{align*}
where $w(\xi)=u(x)$. Therefore, if we set $\widetilde \varphi(\xi):=\|\nabla w(\xi)\|_{\R^d}$, then $\widetilde \varphi(\xi)=\varphi(x)$ and
\begin{align*}
\widetilde \varphi_\varepsilon(\xi)
:= & \sqrt{\varepsilon+\|\nabla_\xi w(\xi)\|_{\R^d}^2}= \varphi_\varepsilon(x),
\end{align*}
with $x=R\xi$. Finally, let us notice that at $\overline x\in \partial O$ we have
\begin{align*}
R^{-1}\nu_{\partial O}(\overline x)
= \nu_{\partial R^{-1}(O)}(\overline \xi)
= & (0,\ldots,0,-1)=(\nabla\psi(\overline \xi'),-1)
\end{align*}
in the system $\xi$, and so $H^\gamma_{\partial R^{-1}(O)}(\overline \xi)=\Delta_{\xi'}\psi(\overline \xi')-\langle \overline \xi,(\nabla \psi(\overline \xi'),-1)\rangle_{\R^d}=\Delta_{\xi'}\psi(\overline \xi')+\overline \xi_d$.

\begin{lem}
\label{lem:graph1}
Let $y\in C^\infty_c(O)$, let $\sigma>0$, let $u=J_\sigma(y)$, let $\overline x\in\partial O$ and let $\overline \xi$ be the corresponding point to $\overline x$ under the rotation coordinates $x=R \xi$ with $R$ as above. Then,
\[
D_{\xi_d}\widetilde \varphi_{\varepsilon}(\overline \xi)= (\widetilde \varphi_{\varepsilon}(\overline \xi))^{-1}
(D_{\xi_d}v(\overline \xi))^{2}\left(\Delta_{\xi'}\psi(\overline \xi')+\overline \xi_{d}\right),
\]
where $\overline \xi,v,\psi$ and $\widetilde \varphi$ are as above.
\end{lem}
\begin{proof}
We replicate the argument of \cite[Lemma 8.2]{Bar}. Let $v(\xi)=u(R^{-1}\xi)$ for every $\xi\in \R^d$. Since $0=u(x)=v(\xi',\psi(\xi'))$ for every $\xi'\in V\subseteq \R^{d-1}$ with $x=R(\xi',\psi(\xi'))$, it follows that
\begin{equation}
D_{\xi_i}v(\xi)+D_{\xi_d}v(\xi)D_{\xi'_i}\psi(\xi')=0
\label{eq:first_deriv1}
\end{equation}
 for every $i\in\{1,\ldots,d-1\}$ and $\xi=(\xi',\psi(\xi'))$ with $\xi'\in V$. Differentiating \eqref{eq:first_deriv1} with respect to $i\in\{1,\ldots,d-1\}$ it follows that
\[
D_{\xi_i\xi_i}^{2}v(\xi)+2D^2_{\xi_i\xi_d}v(\xi)D_{\xi'_i}\psi(\xi')+D_{\xi_d\xi_d}^{2}v(\xi)(D_{\xi'_i}\psi(\xi'))^{2}+D_{\xi_d}v(\xi)D_{\xi'_i\xi'_i}^{2}\psi(\xi')=0
\]
for every $\xi=(\xi',\psi(\xi'))\in (V,\psi(V))$. If we choose $\xi=\overline \xi$ then we get $D_{\xi'_i}\psi(\overline \xi')=0$. As a byproduct, $D_{\xi_i}v(\overline \xi)=0$ and 
\[
D_{\xi_i\xi_i}^{2}v(\overline \xi)+D_{\xi_d}v(\overline \xi)D_{\xi'_i\xi'_i}^{2}\psi(\overline \xi')=0
\]
for every $i\in\{1,\ldots,d-1\}$. By recalling that $\overline L$ can be applied to
$u$ also up the boundary since $u\in C^{2}_0(\overline{O})$ and
recalling that $D_{\xi_i}v(\overline \xi')=0$ for every $i=1,\ldots,d-1\}$, we get
\[
\overline L v(\overline \xi)=\Delta_\xi v(\overline \xi)-\overline \xi_d D_{\xi_d}v(\overline \xi)=D_{\xi_d\xi_d}^{2}v(\overline \xi)-D_{\xi_d}v(\overline \xi)\Delta_{\xi'}\psi(\overline \xi')-\overline \xi_{d}D_{\xi_d}v(\overline \xi).
\]
We notice that $\overline Lv(\overline \xi)=0$ since $\sigma \overline Lv(\overline \xi)=\sigma \overline Lu(\overline x)=u(\overline x)-y(\overline x)=0$, and so 
\[
D_{\xi_d\xi_d}^{2}v(\overline \xi)=D_{\xi_d}v(\overline \xi)(\Delta_{\xi'}\psi(\overline \xi')+\overline \xi_{d}).
\]
Finally, since 
\begin{align*}
\widetilde \varphi_\varepsilon(\overline \xi)D_{\xi_d}\widetilde \varphi_\varepsilon (\overline \xi)=
\sum_{i=1}^dD_{\xi_i}v(\overline \xi)D^2_{\xi_i\xi_d}v(\overline \xi)
= D_{\xi_d}v(\overline \xi)D^2_{\xi_d\xi_d}v(\overline \xi)
\end{align*}
and $u$ (and consequently $v$) is smooth up to the boundary $\partial O$, we conclude that
\begin{align*}
D_{\xi_d}\widetilde \varphi_\varepsilon (\overline \xi)
= & (\widetilde \varphi_\varepsilon(\overline \xi))^{-1}D_{\xi_d}v(\overline \xi)D^2_{\xi_d\xi_d}v(\overline \xi) \\
= & (\widetilde \varphi_\varepsilon(\overline \xi))^{-1}(D_{\xi_d}v(\overline \xi))^2(\Delta_{\xi'}\psi(\overline \xi')+\overline \xi_{d}).
\end{align*}
\end{proof}
Now we state the main result of this section.

\begin{prop}
\label{prop:Diric_21}
Under Hypothesis \ref{hyp:gauss_curv1}, for every $p\in[1,\infty)$ we have 
\begin{align}
\label{eq_fond_gener1}
\int_O\|\nabla J_\sigma y\|^p_{\R^d}d\gamma\leq \int_O\|\nabla y\|_{\R^d}^pd\gamma,
\end{align}
for every $y\in W^{1,q}_0(O,\gamma)$, where $q=p$ if $p>1$ and $q>1$ if $p=1$. 
In particular, $L_{O,p}$ is dissipative in $W^{1,p}_0(O,\gamma)$ for every $p\in(1,\infty)$.
\end{prop}
\begin{proof}
The proof is analogous to that of \cite[Proposition 8.2]{Bar}, hence we skip the details. We assume that $g(t)=t^p$ belongs to $C^2([0,\infty))$ (if $p<2$ then we take a smooth convex approximation), and we consider $y\in C^\infty_c(O)$ and $\sigma>0$. For every $\varepsilon >0$ we set $\psi(x)=g(\varphi(x))$ and $\psi_\varepsilon(x):=g(\varphi_\varepsilon(x))$ for every $x\in O$, where $\varphi$ and $\varphi_\varepsilon$, with $u=J_\sigma(y)$, have been defined in \eqref{def_varphi}. From \eqref{eq:lemma1} and the properties of $g$ we get
\begin{align}
\label{stima_grad_convesso1}
\frac{\psi^2}{\psi_\varepsilon}-\sigma \overline L\psi_\varepsilon\leq 
g'(\varphi_\varepsilon)\left(\|\nabla y\|_{\R^d}-\frac{\varphi^2}{\varphi_\varepsilon}\right)+\frac{\psi^2}{\psi_\varepsilon}
\end{align}
in $O$. Further, for every $\overline x\in \partial O$, if $\overline \xi$ is the corresponding point under the rotation $R^{-1}$ introduced above, Lemma \ref{lem:graph1} gives
\[
\frac{\partial\varphi_{\varepsilon}}{\partial \nu_{\partial O}}(\overline x)
= -D_{\xi_d}\widetilde  \varphi_\varepsilon(\overline \xi)
= - (\widetilde \varphi_\varepsilon(\overline \xi))^{-1}(D_{\xi_d}v(\overline \xi))^2(\Delta_{\xi'}\psi(\overline \xi')+\overline \xi_{d})
=-H'_{\partial R^{-1}(O)}(\overline \xi)\leq 0
\]
for every $\overline x\in\partial O$. This implies that
\begin{align*}
\frac{\partial \psi_\varepsilon}{\partial \nu_{\partial O}}(\overline x)
= -D_{\xi_d}\left(g(\widetilde \varphi_\varepsilon)\right)(\overline \xi)
= -g'(\widetilde \varphi_\varepsilon(\overline\xi))D_{\xi_d}\widetilde \varphi(\overline\xi)\leq 0
\end{align*}
for every $\overline x\in\partial O$. By applying the divergence theorem and taking into account the explicit formula of $\overline L$ we get
\begin{align}
\label{segno_int_bordo1}
\int_{O}\overline L\psi_{\varepsilon}\ d\gamma=\int_{\partial O}\frac{\partial\psi_{\varepsilon}}{\partial \nu_{\partial O}}\theta_d\ dH^{d-1}\leq0,
\end{align}
where $\theta_d$ is the density of the standard Gaussian measure in $\R^{d}$ and $H^{d-1}$
is the $(d-1)$-Hausdorff measure. 
From \eqref{stima_grad_convesso1} and \eqref{segno_int_bordo1} it follows that
\begin{align*}
\int_O \frac{\psi^2}{\psi_\varepsilon}d\gamma
\leq \int_O\left(g'(\varphi_\varepsilon)\left(\|\nabla y\|_{\R^d}-\frac{\varphi^2}{\varphi_\varepsilon}\right)+\frac{\psi^2}{\psi_\varepsilon}\right)d\gamma.
\end{align*}
Letting $\varepsilon$ tend to $0$ we get
\begin{align*}
\int_O g(\|\nabla u\|_{\R^d})d\gamma
\leq & \int_O\left(g'(\varphi)\left(\|\nabla y\|_{\R^d}-\varphi\right)+g(\varphi)\right)d\gamma.
\end{align*}
The convexity of $g$ implies that $g'(u)\left(v-u\right)+g(u)\leq g(v)$ for every $u,v\in[0,\infty)$, which gives \eqref{eq_fond_gener1} for every $y\in C_c^\infty(O)$. From the density of this set in $W^{1,p}_0(O,\gamma)$ and Remark \ref{rem:lower_semicont} we get the thesis.
\end{proof}

\begin{rem}
\label{rem:non_p_q_1}
The fact that the thesis does not holds for $p=q=1$ follows from the fact that, for $r=1$, the arguments in Remark \ref{rem:lower_semicont} do not work.     
\end{rem}

\section{Gradient resolvent contractivity}
\label{sec:res_contr_infin_dim}

In this section we consider a separable Banach space $X$ with nondegenerate centered Gaussian measure $\gamma$ and Cameron-Martin space $H$. We fix an orthonormal basis $\Phi:=\{h_{n}:n\in\N\}$ of $H$ of elements of $QX^*$, that $F_{m}:={\rm span}\{h_{1},\ldots,h_{m}\}$ for every $m\in\N$, that $\pi_{F_m}$ is the projection on $F_m$ and that $X_{m}^{\perp}:={\rm Ker}(\pi_m)$. We can identify $F_{m}$ with $\R^{m}$ by means of the operator $\Pi_m:F_m\rightarrow \R^m$, defined by
\begin{align*}
\Pi_m(y)=(y_1,\ldots,y_m), \quad y\in F_m, \ \ y_i:=\langle y,h_i\rangle_H, \ i=1,\ldots,m.
\end{align*}

We split the proof of the main result into two parts: in the former we consider the case when $O$ is a cylindrical set with respect to the basis $\{h_n:n\in\N\}$, in the latter we consider a generic open set $O\subseteq X$ which satisfies suitable conditions.

\subsection{\label{sub:Cylindrical-case}Cylindrical case}
In this subsection we assume the following additional hypothesis.
\begin{hyp}
\label{hyp:dom_cyl_case}
$O\subseteq X$ is a $C^{2,\alpha}$-regular cylindrical set with respect to the basis $\{h_n:n\in\N\}\subseteq QX^*$. This means that there exist ${\overline m}\in\N$ and a $C^{2,\alpha}$-regular open set $\mathscr O_{\overline m}\subseteq \R^{\overline m}$ such that $O=(\Pi_{{\overline m}}\circ \pi_{\overline m})^{-1}\left(\mathscr O_{{\overline m}}\right)=(\Pi_{\overline m})^{-1}(\mathscr O_{\overline m})\oplus X^\perp_{\overline m}$, and $\mathscr O_{{\overline m}}=g^{-1}((-\infty,0))$ for some function $g\in C^{2,\alpha}(\R^{\overline m})$.
\end{hyp}

\begin{defn}
If we set $G(x):=g((\Pi_{\overline m}\circ\pi_{\overline m})x)$ for every $x\in X$, then it follows that $O=G^{-1}((-\infty,0))$, and the spaces $W^{1,p}(O,\gamma)$ and $W_{0}^{1,p}(O,\gamma)$ introduced in Section \ref{sub:sobolev_spaces} are well-defined. Finally, we notice that for every $n\geq {\overline m}$ we have $\mathscr O_n=\mathscr O_{\overline m}\times \R^{n-{\overline m}}$.
\end{defn}

\begin{defn}
For every $k\in\N$ we denote by $\mathcal{F}C^{k}_{b,\Phi}(X)$ the set of cylindrical $C^{k}_b$ functions with respect to the orthonormal basis $\Phi$, i.e., a function $y\in\mathcal{F}C^{k}_{b,\Phi}(X)$ if there exists $n\in \N$ and $v\in C_b^k(\R^n)$ such that $y(x)=v((\Pi_n\circ\pi_n)(x))$ for every $x\in X$.

We denote by $\mathcal{F}C_{b,0}^{k}(O)$ the subset of $\mathcal{F}C_{b,\Phi} ^k(X)$ of functions which vanish out of an open set $A$ with positive distance from $O^c$, i.e., a function $y$ on X belongs to $\mathcal{F}C_{b,0}^{k}(O)$ if $y\in\mathcal {F}C^k_{b,\Phi}(X)$ and there exists an open set $A$, with $d(A,O^c)>0,$ such that $y=0$ on $A^c$. In particular, from the definition of $O$ it follows that every $y\in\mathcal FC_{b,0}^k(O)$ has the form $y=v(\Pi_n\circ\pi_n)$ for some $n\geq \overline m$ and $v\in C_b^k(\R^n)$.

\end{defn}

\begin{lem}
\label{lem:dens_cyl_func_cyl_dom}
If $O$ is cylindrical with respect to the basis $\Phi$, then $W^{1,p}_0(O,\gamma)$ is the closure of $\mathcal{F}C_{b,0}^1(O)$ with respect to the norm of $W^{1,p}(O,\gamma)$.
\end{lem}
\begin{proof}
At first we show that if $y\in \mathcal {F}C_{b,0}^k(O)$ and $v\in C_b^k(\R^n)$ satisfies $y(x)=v((\Pi_n\circ\pi_n)x)$ for every $x\in X$, then $v$ has support in an open set $\mathcal A_n\subseteq \mathcal O_n$ with positive distance from $\mathcal O_n^c$.

Let $A\subseteq O$ be an open set with positive distance from $O^c$ such that $y$ vanishes out of $A$. We are going to prove that the above assertion is true with $\mathcal A_n:=\{(\Pi_n\circ\pi_n)x:x\in A\}$. 

\noindent At first we show that $\mathcal A_n$ is an open set. Let $x\in A$, $x_n:=(\Pi_n\circ\pi_n)x\in \mathcal A_n$ and let $\varepsilon>0$ be such that $B(x,\varepsilon)=\{y\in X:\|y-x\|_X<\varepsilon\}\subseteq A$. For every $y_n\in \R^n$ such that $\|y_n-x_n\|_{\R^n}<c_H^{-1}\varepsilon$ (where $c_H$ is the smallest constant in the continuous embedding $H\subseteq X$, see Subsection \ref{sub:Fundamentals-about-Wiener} and formula \eqref{cost_C_H}), we set $y=(\Pi_n\circ \pi_n)^{-1}y_n+(x-(\Pi_n\circ\pi_n)^{-1}x_n)\in X$ and we get
\begin{align*}
\|y-x\|_X=\|(\Pi_n\circ \pi_n)^{-1}(y_n-x_n)\|_X\leq c_H\|\Pi^{-1}_n(y_n-x_n)\|_H=c_H\|y_n-x_n\|_{\R^n}< \varepsilon.    
\end{align*}
Hence, $y\in A$ and so $y_n=(\Pi_n\circ\pi_n)y\in \mathcal A_n$. This means that the ball centered at $x_n$ with radius $c_H^{-1}\varepsilon$ in $\R^n$ is contained in $\mathcal A_n$. 

\noindent It remains to show that $v$ vanishes on $\mathcal A_n^c\subseteq \R^n$. To this aim, we fix $\xi=(\xi_1,\ldots,\xi_n)\in \mathcal A_n^c$ and set $w:=\sum_{i=1}^n\xi_i h_i\in F_{n}$. If there exists $z\in X^\perp_n$ such that $x=w+z\in A$, then $\xi=(\Pi_n\circ\pi_n)x\in \mathcal A_n$, which gives a contradiction. This implies that $\{w\}\oplus X^\perp_n \subseteq A^c$, and so $v(\xi)=0$.

Now we are able to prove the statement. Let $f\in {\mathcal H}_{b,0}^1(O)$, and for every $n\geq {\overline m}$ let us set $f_n:=\mathbb E[\overline f|F_n]$, where
\begin{align*}
f_n(x):=\int_X\overline f(\pi_n(x)+(Id-\pi_n)(y))d\gamma(y), \quad x\in X.
\end{align*}
Let $A\subseteq O$ be an open set with positive distance from $O^c$ such that ${\rm supp}(f)\subseteq A$. We set $\widetilde A:=\{\pi_{\overline m}(x)+y:x\in A, \ y\in X_{\overline m}^\perp\}=\{x\in X:\pi_{\overline m}(x)=\pi_{\overline m}(z) \textrm{ for some } z\in A\}\subseteq X$. We claim that $A\subseteq \widetilde A\subseteq O$ is open and has positive distance from $O^c$.

\noindent The inclusion $A\subseteq \widetilde A$ is trivial. To show that $\widetilde A\subseteq O$, we notice that if $\pi_{\overline m}(x)+y\in O^c$ for some $x\in A$ and $y\in X_{\overline m}^\perp$, then from the definition of $O$ we get $\pi_{\overline m}(x)+z\in O^c$ for every $z\in X_{\overline m}^\perp$. In particular, the choice of $z=x-\pi_{\overline m}(x)$ implies that $x\in O^c$, which contradicts the fact that $x\in A$. 

\noindent Now we prove that $\widetilde A$ is open. let us fix $\widetilde x=\pi_{\overline m}(x)+y\in \widetilde A$, where $x\in A$ and $y\in X_{\overline m}^\perp$. Since $A$ is open, there exists $\varepsilon>0$ such that $B(x,\varepsilon)=\{z\in X:\|z-x\|_X<\varepsilon\}\subseteq A$. For every $\widetilde y\in B(\widetilde x,\varepsilon)$, from the decomposition $X=F_{\overline m}\oplus X_{\overline m}^\perp$ we get 
\begin{align*}
\varepsilon>\|\widetilde x-\widetilde y\|_X\geq \|\pi_{\overline m}(\widetilde x-\widetilde y)\|_X=\|\pi_{\overline m}(x)-\pi_{\overline m}(\widetilde y)\|_X.
\end{align*}
If we consider $y=\pi_{\overline m}(\widetilde y)+(x-\pi_{\overline m}(x))$, it follows that $\pi_{\overline m}(y)=\pi_{\overline m}(\widetilde y)$ and $\|x-y\|_X=\|\pi_{\overline m}(x)-\pi_{\overline m}(\widetilde y)\|_X<\varepsilon$, which means that $y\in A$.

\noindent It remains to show that $\widetilde A$ has positive distance from $O^c$. For this purpose, for every $\widetilde x\in \widetilde A$ and every $z\in O^c$ we get
\begin{align*}
\|\widetilde x-z\|_X
\geq \|\pi_{\overline m}(x)-\pi_{\overline m}(z)\|_X = \|x-(\pi_{\overline m}(z)+x-\pi_{\overline m}(x))\|_X,
\end{align*}
where $x\in A$ satisfies $\pi_{\overline m}(x)=\pi_{\overline m}(\widetilde x)$. Since $z\in O^c$, the definition of $O$ gives $\pi_{\overline m}(z)+x-\pi_{\overline m}(x)\in O^c$, and so $ \|x-(\pi_{\overline m}(z)+x-\pi_{\overline m}(x))\|_X\geq {\rm dist}(A,O^c)$. The arbitrariness of $\widetilde x\in \widetilde A$ and $z\in O^c$ imply that ${\rm dist}(\widetilde A,O^c)\geq {\rm dist}(A,O^c)$. The claim is proved.

Let us consider $x\in \widetilde A^c$. From the claim, it follows that $\pi_{\overline m}(x)+y\in A^c$ for every $y\in X_{\overline m}^\perp$. Indeed, if there exists $y\in X_{\overline m}^\perp$ such that $\pi_{\overline m}(x)+y\in A$, we infer that $\pi_{\overline m}(x)=\pi_{\overline m}(\pi_{\overline m}(x)+y)$, which means that $x\in \widetilde A$. Hence, for every $x\in \widetilde A^c$ and $y\in X$ it follows that
\begin{align*}
\pi_n(x)+(Id-\pi_n)(y)=\pi_{\overline m}(x)+(\pi_{n}(x)-\pi_{\overline m}(x))+(Id-\pi_{n})(y)\in A^c,
\end{align*}
which implies that $f_n=0$ on $\widetilde A^c$. In particular, ${\rm supp}(f_n)\subseteq \widetilde A$ and $f_n(x)=f_n(\pi_nx)= v_n((\Pi_n\circ\pi_n)x)$ for every $x\in X$, where $v_n=f_n\circ(\Pi_n\circ \pi_n)^{-1}\in C^1_{b,0}(\mathcal O_n)$ from the first part of the proof. Finally, $(f_n)_{n\in\N}$ converges to $\overline f$ in $W^{1,p}(X,\gamma)$ as $n\to \infty$ from \cite[Corollary 3.5.2]{Bog}, which implies that $({f_n}_{|O})_{n\in\N}$ tends to $f$ in $W^{1,p}(O,\gamma)$ as $n\to \infty$.
\end{proof}

The next proposition shows that the result in Proposition \ref{prop:Diric_21} can be easily generalized in infinite dimension when $O$ is a cylindrical domain.
\begin{prop}
\label{prop:Dirichlet_cylindrical} Let $O$ satisfy Hypothesis \ref{hyp:dom_cyl_case}, and assume that $H^{\gamma^m}_{\partial \mathscr O_{{\overline m}}}(\xi)\geq0$ for every $\xi\in \partial \mathscr O_{\overline m}$. If $p\in[1,\infty)$, then for every $\sigma>0$ we get
\begin{align*}
\int_{O}\|\nabla_{H}J^O_{\sigma}(y)\|^p_{H} d\gamma\leq\int_{O}\|\nabla_{H}y\|_{H}^p d\gamma,
\quad y\in W_{0}^{1,q}(O,\gamma),
\end{align*}
where $q=p$ if $p>1$ and $q>1$ if $p=1$.
 \end{prop}
\begin{proof}
For every $n\geq {\overline m}$ and every $ \xi\in \partial \mathscr O_n$ we have $ \xi=(\overline \xi_0, \xi_1,\ldots, \xi_{n-{\overline m}})$, with $\overline \xi_0\in \partial \mathscr O_{\overline m}$ and $(\xi_1,\ldots, \xi_{n-{\overline m}})\in \R^{n-{\overline m}}$. Further, the outer normal $\nu_n(\xi)$ to $\mathcal O_n$ at $\xi\in\partial \mathcal O_n$ satisfies $\nu_n( \xi)=(\nu_{\overline m}(\overline \xi_0),0,\ldots0)$, where $\nu_{\overline m}(\overline \xi_0)$ is the outer normal to $\mathcal O_{\overline m}$ at $\overline \xi_0\in\partial \mathcal O_{\overline m}$. Hence, we get 
\begin{align*}
H^{\gamma^n}_{\partial \mathscr O_n}(\xi) =H^{\gamma^m}_{\partial \mathscr O_{\overline m}}(\overline \xi_0)\geq0, \quad  \xi\in \partial \mathscr O_n.
\end{align*}

Let us fix $\sigma>0$. On $L^{2}(\mathscr O_{n},\gamma^{n})$, where $\gamma^n$ is the standard Gaussian measure on $\R^n$, we consider the Ornstein-Uhlenbeck operator $L_{\mathscr O_n}$ with homogeneous  Dirichlet boundary conditions with respect to $\mathscr O_{n}$ and the bounded operator $\mathcal J_{\sigma,n}=(I-\sigma L_{\mathscr O_n})^{-1}$. Let $y\in \mathscr {F}C_{b,0}^1(O)$ and let $v\in C_{b,0}^1(\mathscr O_n)$ be such that $y(x)=v((\Pi_n\circ \pi_n)(x))$ for every $x\in X$. We claim that
\[
(\mathcal J_{\sigma,n}v)(\Pi_n\circ\pi_{{n}})=J_{\sigma}^Oy, \quad \gamma \textrm{-a.e. in} \ O. 
\]
Indeed, for every $f\in\mathcal{F}C_{b,0}^1(O)$ with $f=\varphi(\Pi_j\circ\pi_j)$, where $j\in\N$ and $\varphi\in C^1_b(\mathbb R^j)$, if $\overline m\leq j\leq n$ then we get
\begin{align*}
\int_O\langle \nabla_H\mathcal J_{\sigma,n}v(\Pi_n\circ\pi_n),\nabla_Hf\rangle_H d\gamma
= & \int_{\mathcal O_n}\langle \nabla \mathcal J_{\sigma,n}v,\nabla \varphi\rangle_{\R^n}d\gamma^n 
= \int_{\mathcal O_n}L_{\mathcal O_n}(\mathcal J_{\sigma,n}v)\cdot \varphi d\gamma^n \\
= & \int_O (L_{\mathcal O_n}(\mathcal J_{\sigma,n}v))(\Pi_n\circ \pi_n)\cdot f d\gamma,
\end{align*}
and if $j>n$ then it follows that
\begin{align*}
\int_O\!\langle \nabla_H\mathcal J_{\sigma,n}v(\Pi_n\circ\pi_n),\nabla_Hf\rangle_H d\gamma
= & \int_{\R^{j-n}}\int_{\mathcal O_n}\!\langle \nabla \mathcal J_{\sigma,n}v(\xi),\nabla \varphi(\xi+\eta)\rangle_{\R^n}d\gamma^n(\xi) d\gamma^{j-n}(d\eta) \\
= & \int_{\R^{j-n}}\int_{\mathcal O_n}L_{\mathcal O_n}(\mathcal J_{\sigma,n}v)(\xi)\cdot \varphi(\xi+\eta) d\gamma^n(\xi)d\gamma^{j-n}(\eta) \\
= & \int_O (L_{\mathcal O_n}(\mathcal J_{\sigma,n}v))(\Pi_n\circ \pi_n)\cdot f d\gamma.
\end{align*}
From the definition of $L_O$ and Lemma \ref{lem:dens_cyl_func_cyl_dom} we conclude that $(\mathcal J_{\sigma,n}v)(\Pi_n\circ \pi_n)\in D(L_O)$ and 
\begin{align*}
L_O((\mathcal J_{\sigma,n}v)(\Pi_n\circ \pi_n))=(L_{\mathcal O_n}(\mathcal J_{\sigma,n}v))(\Pi_n\circ \pi_n).
\end{align*}
Moreover, for $\gamma$-a.e. $x\in X$ we get
\begin{align*}
\sigma L_O((\mathcal J_{\sigma,n}v)(\Pi_n\circ \pi_n))(x)
= & \sigma L_{\mathcal O_n}((\mathcal J_{\sigma,n}v))(\Pi_n\circ \pi_n)(x) \\
= & (\mathcal J_{\sigma,n}v)(\Pi_n\circ \pi_n)(x)-v(\Pi_n\circ\pi_n)(x).
\end{align*}
Combined with the fact that $y=v(\Pi_n\circ\pi_n)$, this gives the claim. In particular, we get
$\|\nabla_HJ_\sigma y\|_{L^p(O,\gamma,H)}=\|\nabla \mathcal J_{\sigma,n}v\|_{L^p(\mathcal O_n,\gamma^n,\R^n)}$ for every $y\in \mathcal FC_{b,0}^1(O)$.

Since $\|\nabla_Hy\|_{L^p(O,\gamma,H)}=\|\nabla v\|_{L^p(\mathcal O_n,\gamma^n,\R^n)}$, from Proposition \ref{prop:Diric_21} it follows that
\[
\int_{O}\|\nabla_{H}J_{\sigma}y\|_{H}^pd\gamma
=\!\int_{\mathscr O_{n}}\|\nabla \mathcal J_{\sigma,n}v\|_{\R^n}^pd\gamma^{n}
\leq\!\int_{\mathscr O_{n}}\|\nabla v\|_{\R^n}^p d\gamma^{n}
=\!\int_{O}\|\nabla_{H}y\|_H^pd\gamma,
\]
and we get the thesis when $y\in \mathscr{F}C_{b,0}^{1}(O)$. From Remark \ref{rem:lower_semicont} and Lemma \ref{lem:dens_cyl_func_cyl_dom} we conclude.
%
\end{proof}

\subsection{\label{sub:Generalization-to-non-cylindrica}Generalization to the non-cylindrical
case}
For every $f\in \mathcal{F}C _b^\infty(X)$ of the form $f=v(\Pi_n\circ\pi_n)$, with $v\in C_b^\infty(\R^n)$, we have
\begin{align*}
Lf=\sum_{i=1}^n\partial^2_{h_ih_i}f-\sum_{i=1}^n\partial_{h_i}f\hat h_i.
\end{align*}
Let us notice that for these functions we have
\begin{align*}
Lf(x)=L_n v(\xi):=\Delta v(\xi)-\langle \xi, \nabla v(\xi)\rangle_{\R^n}, \quad \xi=(\Pi_n\circ\pi_n)(x), \ x\in X,
\end{align*}
where $\Delta$ is the Laplace operator in $\R^n$.

Assume that Hypothesis \ref{claim:regularity} are satisfied, and recall that $\partial O=G^{-1}(0)$ and that at every point  $x\in \partial O$ the outer $H$-normal is $\nabla_HG(x)/\|\nabla_HG(x)\|_H$. For every $n\in\N$ we consider the cylindrical function $G_n:=G\circ \pi_n$ (see \cite[Section 4]{Dap}). From the definition, $G_n$ is a $H$-differentiable cylindrical function and $\nabla_HG_n$ is everywhere $H$-differentiable. Hence, there exists $g_n\in C^2(\R^n)$ such that $G_n(x)=g_n((\Pi_n\circ\pi_n)(x))$ for every $x\in X$. For every $n\in\N$ we introduce the sets $O_n:=G_n^{-1}((-\infty,0))$ and $\mathcal O_n:=g_n^{-1}((-\infty,0))$. Therefore, $\partial O_n=G_n^{-1}(0)$, $\mathcal O_n=g_n^{-1}(\{0\})$ and $O_n=(\Pi_n\circ \pi_n)^{-1}(\mathscr O_n)$ for every $n\in\N$. For every $x\in \partial O_n$ with $\nabla_HG_n(x)\neq 0$ we set
\begin{align*}
H_{\partial O_n}^\gamma(x):=\frac{LG_n(x)}{\|\nabla_HG_n(x)\|_H}-\frac{\langle D^2_HG_n(x)\nabla_HG_n(x), \nabla_HG_n(x)\rangle}{\|\nabla_HG_n(x)\|_H^3}.
\end{align*}

Let us assume the following conditions on $O_n$, $\mathscr O_n$ and $H_{\partial O_n}^\gamma$.
\begin{hyp}
\label{claim:Gaussian-curvature}We suppose that there exists $n_0\in \N$ such that $\|\pi_{n_0}(\nabla_{H}G)\|_H\neq 0$ on $\partial O$, for every $n\geq n_0$ the set $\mathscr O_n$ is a  $C^{2,\alpha}$-regular open set in $\R^n$ and $H_{\partial O_n}^\gamma(x)\geq0$ for every $x\in\partial O_{n}$.
\end{hyp}

\begin{rem}
\label{rmk:corvatura_n}
Assume that Hypothesis \ref{claim:Gaussian-curvature} is satisfied. At first, we notice that, for every $n\geq n_0$, the assumption $\|\pi_{n_0}(\nabla_{H}G)(x)\|_H\neq 0$ for every $x\in \partial O$ implies that, for every $x\in\partial O_n$, we get $\|\nabla_H G_n(x)\|_H\geq \|\pi_{n_0}(\nabla_HG_n(x))\|_H=\|\pi_{n_0}(\nabla_HG(\pi_n(x)))\|_H\neq 0$, since $x\in \partial O_n$ if and only if $\pi_n(x)\in \partial O$.  

Hence, $\nu_{\partial \mathscr O_n}=\nabla g_n/\|\nabla g_n\|_{\R^n}$ for every $n\geq n_0$,  and for every $x\in \partial O_n$ we have
\begin{align*}
\frac{\langle D_H^2 G_n(x)\nabla_H G_n(x),\nabla_H G_n(x)\rangle_H}{\|\nabla_H G_n\|_{H}^3}
=\frac{\langle D^2 g_n(\xi)\nabla g_n(\xi),\nabla g_n(\xi)\rangle_{\R^n}}{\|\nabla g_n\|_{\R^n}^3},
\end{align*}
with $\xi=(\Pi_n\circ\pi_n)(x)$. This implies that 
\begin{align*}
H_{\partial \mathscr O_n}^{\gamma^n}(\xi)
= &  \frac{\Delta_n g_n(\xi)}{\|\nabla g_n(\xi)\|_{\R^n}}-\frac{\langle D^2 g_n(\xi)\nabla g_n(\xi),\nabla g_n(\xi)\rangle_{\R^n}}{\|\nabla g_n\|_{\R^n}^3}-\frac{\langle \xi,\nabla g_n(\xi)\rangle_{\R^n}}{\|\nabla g_n(\xi)\|_{\R^n}}  \\
= & \frac{L_ng_n(\xi)}{\|\nabla g_n(\xi)\|_{\R^n}}-\frac{\langle D^2 g_n(\xi)\nabla g_n(\xi),\nabla g_n(\xi)\rangle_{\R^n}}{\|\nabla g_n\|_{\R^n}^3}  \\
= & \frac{LG_n(x)}{\|\nabla_HG_n\|_H}-\frac{\langle D_H^2 G_n(x)\nabla_H G_n(x),\nabla_H G_n(x)\rangle_H}{\|\nabla_H G_n\|_{H}^3} \\
= & H_{\partial O_n}^\gamma(x)\geq0,
\end{align*}
for every $x\in \partial O_n$ and $\xi:=(\Pi_n\circ\pi_n)(x)\in\partial\mathscr O_n$.
\end{rem}

\begin{thm}
\label{thmgeneral_case} Let $O\subseteq X$ be an open set satisfying Hypotheses \ref{claim:regularity} and  \ref{claim:Gaussian-curvature}. Then, for every $p\in[1,\infty)$ and $\sigma>0$, it follows that 
\begin{align*}
\int_{O}\|\nabla_{H}J_{\sigma}(y)\|^p_{H}\ d\gamma\leq\int_{O}\|\nabla_{H}y\|^p_{H}\ d\gamma,
\end{align*}
for every $y\in W_{0}^{1,q}(O,\gamma)$, where $q=p$ if $p>1$ and $q>1$ if $p=1$. In particular, the operator $L_{O,p}$ is dissipative in $W^{1,p}_0(O,\gamma)$ for every $p\in(1,\infty)$.
\end{thm}

\begin{rem}
For the case $p=q=1$ see Remark \ref{rem:non_p_q_1}.    
\end{rem}

\begin{proof}
Let $p$ and $q$ be as in the statement.

Let $y\in \mathcal H_{b,0}^1(O)$ and for every $n\in\N$ let us set $y_n:=y\circ \pi_n$. From Lemma \ref{lem:H_O_n_da_H_O} we deduce that $y_n\in \mathcal H_{b,0}^1(O_n)$, and from Remark \ref{rmk:corvatura_n} it follows that $H_{\partial {\mathcal O}_n}^\gamma\geq0$ on $\partial\mathscr O_n$, which implies that the assumptions of Lemma \ref{prop:Dirichlet_cylindrical} are satisfied. Let $\sigma>0$. Since $y_n\in \mathcal H_{b,0}^1(O_n)\subseteq W^{1,q}_0(O_n,\gamma)$ we get 
\begin{align*}
\int_{O_{n}}\|\nabla_{H}J_{\sigma,n}(y_n)\|^p_{H}\ d\gamma\leq\int_{O_{n}}\|\nabla_{H}y_n\|^p_{H}\ d\gamma,
\end{align*}
for every $n\geq n_0$, where $L_{O_n}$ is the Ornstein-Uhlenbeck operator with homogeneous Dirichlet boundary conditions on $O_{n}$, and $J_{\sigma,n}=({\rm Id}-\sigma L_{O_n})^{-1}$. 

Since $y\in W^{1,q}_0(O,\gamma)$ it follows that $\overline y\in W^{1,q}(X,\gamma)$ and $\overline {\nabla_Hy}(x)=\nabla _H\overline y(x)$ for $\gamma$-a.e. $x\in X$. Analogously, $\overline y_n\in W^{1,q}(X,\gamma)$ and $\overline {\nabla_Hy_n}(x)=\nabla _H\overline {y_n}(x)$ for $\gamma$-a.e. $x\in X$, for every $n\in\N$. Since $\overline y\in \mathcal H^1(X)$ is bounded and $\overline {y_n}=\overline y\circ \pi_n$ for every $n\in\N$, from Lemma \ref{lem:convergence_projec_1} we infer that $ (\overline{y_n})_{n\in\N}$ converges to $\overline y$ in $W^{1,q}(X,\gamma)$ as $n$ tends to infinity, and so
\begin{align*}
\int_O\|\nabla_Hy\|^p_Hd\gamma
=& \int_X\|\overline {\nabla _Hy}\|^p_H d\gamma
= \int_X\| \nabla_H \overline y\|^p_Hd\gamma
=\lim_{n\rightarrow \infty}
\int_X\| \nabla_H \overline {y_n}\|^p_Hd\gamma \\
= &\lim_{n\rightarrow \infty}
\int_X\| \overline {\nabla_H {y_n}}\|^p_Hd\gamma
=  \lim_{n\rightarrow \infty}
\int_{O_n}\|\nabla_H y_n\|^p_Hd\gamma.
\end{align*}
If we set $u:=J_\sigma^O(y)$ and $u_n:=J_{\sigma,n}^{O_n}(y_n)$ for every $n\geq n_0$, then from Proposition \ref{prop:Dirichlet_cylindrical} applied to $y_n$ we deduce that
\begin{align}
\label{stima_contr_1}
\liminf_{n\rightarrow\infty}
\int_{O_n}\|\nabla_Hu_n\|_H^pd\gamma
\leq \liminf_{n\rightarrow\infty}
\int_{O_n}\|\nabla_H y_n\|_H^pd\gamma
= \int_O\|\nabla_Hy\|_H^pd\gamma.
\end{align}
From Proposition \ref{prop:conv_ris_n} we infer that $(\overline u_n)_{n\in\N}$ converges to $\overline u$ in $W^{1,2}(X,\gamma)$ as $n$ goes to infinity, and Remarks \ref{rem:grad_est_nulla} and \ref{rem:lower_semicont} (with $O=X$ and $r=2$) imply that
\begin{align}
\label{stima_contr_2}
\int_O\|\nabla_Hu\|_H^pd\gamma
= & \int_X\|\nabla_H \overline u\|_H^pd\gamma
\leq \liminf_{n\rightarrow\infty}\int_X\|\nabla _H\overline u_n\|_H^pd\gamma
= \liminf_{n\rightarrow\infty}\int_{O_n}\|\nabla _Hu_n\|_H^pd\gamma.
\end{align}
From \eqref{stima_contr_1} and \eqref{stima_contr_2} we conclude that
\begin{align*}
\int_O\|\nabla_HJ_\sigma(y)\|^p_Hd\gamma
\leq  \int_O\|\nabla_Hy\|^p_Hd\gamma,
\end{align*}
for every $y\in {\mathscr H}_{b,0}^1(O)$. For the general case $y\in W_{0}^{1,q}(O,\gamma)$,
we recall that ${\mathscr H}_{b,0}^1(O)$ is dense in $W_{0}^{1,q}(O)$ by Lemma \ref{lem:dens_H_funct_sob_dir}, and the thesis follows from the boundedness of $J^O_{\sigma}$ on $L^q(O,\gamma)$ and Remark \ref{rem:lower_semicont}.
\end{proof}

\section{Examples}
\label{sec:examples}

In this last section we provide three examples to which our results apply. Such situations have been already considered in \cite[Section 5]{Cel} and \cite[Section 5]{Dap}.

\subsection{Epigraphs}

Let $(X,\gamma)$ be an abstract Wiener space. We fix an orthonormal basis $\{h_{n}:n\in\N\}$ of $H$ in $Q(X^*)$, i.e., $h_n=Q(x^*_n)$ with $x^*_n\in X^*$ for every $n\in\N$. We  define a function $G$ such that $O=G^{-1}((-\infty,0))$ is the epigraph of a function.

\begin{claim}
\label{claim:epigrap}Let $\Phi:X\rightarrow\R$ be a continuous function such that
\begin{enumerate}
\item $\Phi\in \mathcal H^{1}(X)$ (hence $\nabla_{H}\Phi$ exists everywhere);
\item $\partial_{h_{1}}(\Phi)(x)=0$ for every $x\in X$;
\item $\nabla_{H}\Phi$ is everywhere $H$-differentiable, with $\|D_{H}^{2}\Phi\|_{HS}$
uniformly bounded; 
\item $L\Phi$ is uniformly bounded; 
\item the function $f_n:=\Phi\circ(\Pi_n\circ \pi_{n})^{-1}\in C_{\rm loc}^{2,\alpha}(\R^n)$ for every
$n\in\N$, for some $\alpha>0$;
\item there exist $C,C_{1},C_{2},C_{3}>0$, with $C-C_1-C_2-C_3\geq0$, such that $\Phi(x)\geq C$, $\|D^2_H\Phi(x)\|_{HS}\leq C_1$, $\langle \nabla_{H}\Phi(h),h\rangle_H\leq C_{2}$  and $\|D^2_H\Phi(x)\|_{\mathcal L(H,H)}\leq C_3$, for every $x\in X$ and every $h\in H$.
\end{enumerate}
\end{claim}
Hypothesis \ref{claim:epigrap}$(2)$ implies that $\Phi(x)=\Phi(x-\pi_{1}(x))$ for every $x\in X$. We set
\begin{align*}
G(x)=\widehat{h}_{1}(x)+\Phi(x), \quad x\in X,
\end{align*}
and $O:=G^{-1}((-\infty,0))=\{x\in X:\widehat {h}_1(x)\leq -\Phi(x)\}$. 
We have
\begin{align*}
\nabla_{H}G(x)=h_{1}+\nabla_{H}\Phi(x-\pi_{1}(x)), \quad x\in X,
\end{align*}
and
\begin{align*}
D_{H}^{2}G(x)=D_{H}^{2}\Phi(x-\pi_{1}(x)), \quad x\in X.
\end{align*}
By Hypothesis \ref{claim:epigrap}$(1)-(4)$ the function $G$ satisfies Hypothesis \ref{claim:regularity}
(see \cite[Subsection 5.1]{add2}), and for every $n\in\N$ we set
\begin{align*}
G_{n}(x)=\widehat{h}_{1}(x)+\Phi(\pi_{n}(x)-\pi_{1}(x)), \quad x\in X.
\end{align*}
Hence, the function $g_n:=G_n((\Pi_n\circ\pi_n)^{-1})$ belongs to $C^{2,\alpha}(\R^n)$ for every $n\in\N$.
Let us set $\overline x_{n}=\pi_{n}(x)-\pi_{1}(x)\in H$ for every $x\in X$, so $G_{n}(x)=\hat{h}_{1}(x)-\Phi(\overline x_n)$ and 
\begin{align*}
& \nabla_{H}G_{n}(x)=h_{1}+\pi_{n}\nabla_{H}\Phi(\overline x_n), \\
& LG_{n}(x)=\sum_{i=2}^{n}\langle D_{H}^{2}\Phi (\overline x_n) h_{i},h_{i}\rangle_H
-\hat{h}_{1}(x)-\langle\nabla_{H}\Phi(\overline x_n),\overline x_n\rangle _{H},
\end{align*}
This implies that the first part of Hypothesis \ref{claim:Gaussian-curvature} is verified with $n_0=1$, since $\pi_1(\nabla_{H}\Phi(\overline x_n))=0$. Further, on $G^{-1}(0)=\{x\in X:\widehat h_1(x)=-\Phi(x)\}$, for every $n\in\N$ we get 
\begin{align*}
&  H_{\partial O_n}^\gamma(x) 
= \frac{LG_{n}(x)}{\|\nabla_{H}G_{n}(x)\|_H}-\frac{\langle D_{H}^{2}G_{n}(x)\nabla_{H}G_{n}(x),\nabla_{H}G_{n}(x)\rangle_H}{\|\nabla_{H}G_{n}(x)\|_H^{3}} \\
=& \frac{\Phi(x)+\sum_{i=2}^{n}\langle D_{H}^{2}\Phi(\overline x_n)h_{i},h_{i}\rangle_H-\langle\nabla_{H}\Phi(\overline x_n),\overline x_n\rangle _{H}}{\|\nabla_{H}G_{n}(x)\|_H} 
-\frac{\langle D_{H}^{2}\Phi(\overline x_n)\nabla_{H}\Phi(\overline x_n),\nabla_{H}\Phi(\overline x_n)\rangle_H}{\|\nabla_{H}G_{n}(x)\|_H^3} \\
\geq & \frac{C-C_1-C_2}{\|\nabla_{H}G_{n}(x)\|_H}-
\frac{C_3\|\nabla_H\Phi(\overline x_n)\|_H^2}{\|\nabla_{H}G_{n}(x)\|_H^3} 
\geq   \frac{C-C_1-C_2-C_3}{\|\nabla_{H}G_{n}(x)\|_H}\geq 0, \qquad x\in \partial O.
\end{align*}
Hence, $G$ satisfies Hypothesis \ref{claim:Gaussian-curvature} with $n_0=1$. In particular, if $\Phi\equiv C\geq0$ everywhere, the above conditions are verified, which means that the open half-spaces $\{\widehat{h_{1}}<-C\}$ with $C\geq0$
fulfill our assumptions.

\medskip{}

\subsection{Brownian motion and Brownian bridge}


\subsubsection{Brownian motion starting from 0}
\label{ex:br_mot}
We recall the definition of Brownian motion (see \cite[Section 2.3]{Bog}).
We consider the classical Wiener space $(X,\gamma^{W})$ where $X=L^{2}[0,1]$, $\gamma^W$ is the Wiener measure on $L^2[0,1]$, and the Cameron-Martin space $H$ is the set of absolutely continuous functions $f$ on $[0,1]$ such that $f'\in L^{2}[0,1]$ and $f(0)=0$. 

For every $f_{1},f_{2}\in H$ the inner product in $H$ is defined as
\begin{align*}
\left\langle f_{1},f_{2}\right\rangle _{H}=\int_{0}^{1}f_{1}'(x)f_{2}'(x)\ dx.
\end{align*}
We consider the orthonormal basis $\{e_n:n\in\N\}$ of $L^{2}[0,1]$ defined by
\begin{align*}
e_{n}(x):=\sqrt{2}\sin\Bigl(\frac{x}{\sqrt{\lambda_{n}}}\Bigr)=\sqrt{2}\sin\Bigl(\left(n-\frac{1}{2}\right)\pi x\Bigr), \quad x\in[0,1], \ n\in\N,
\end{align*}
where 
\begin{align*}
\lambda_{n}=\frac{1}{\pi^{2}\left(n-\frac{1}{2}\right)^{2}}, \quad n\in\N.
\end{align*}
The system $\{h_{n}=\sqrt{\lambda_{n}}e_{n}:n\in\N\}$ is an orthonormal basis of $H$ and for every $m\in\N$ we denote by $\pi_m:X\rightarrow F_m$ the projection on $F_m:={\rm span}\{h_1,\ldots,h_m\}$. For every $h\in H$ we have
\begin{equation}
\| h\|_{C([0,1])}\leq\int_{0}^{1}|h'(t)|\ dt\leq\left(\int_{0}^{1}|h'(t)|^{2}\ dt\right)^{1/2}=\|h\|_{H}.\label{eq:Bridge_ineq-1}
\end{equation}

Let $g\in C_b^{2,\alpha}(\R)$, $\alpha\in(0,1)$, satisfying for some $c>0$, $\alpha_{1},\alpha_{2},\beta_{1},\beta_{2}\in\R$ the inequalities 
\begin{equation}
|g'(\xi)|\geq c, \quad \alpha_{1}g(\xi)+\beta_{1}\leq\xi g'(\xi)\leq\alpha_{2}g(\xi)+\beta_{2}\label{eq:asymptote-1}
\end{equation}
for every $\xi\in\R$. The above assumptions are satisfied, for instance, by the function $g=p/q$, where $q$ is a positive polynomial of degree $m\in\N$ and $p$ is a polynomial of degree $m+1$
such that $g'(\xi)\neq0$ for all $\xi\in\R$.


\begin{rem}
We have 
\begin{equation}
\sum_{n=1}^{\infty}\left(n-\frac{1}{2}\right)^{-2}=\frac{\pi^{2}}{2}.\label{eq:basel_modified}
\end{equation}
This is a consequence of the well known Basel problem 
\[
\sum_{n=1}^{\infty}n^{-2}=\frac{\pi^{2}}{6},
\]
by using the fact that 
\[
\sum_{n=1}^{\infty}\left(n-\frac{1}{2}\right)^{-2}=4\sum_{i=1}^{\infty}\left(2n-1\right)^{-2}
\]
and 
\[
\sum_{n=1}^{\infty}n^{-2}=\sum_{n=1}^{\infty}\left(2n-1\right)^{-2}+\sum_{n=1}^{\infty}\left(2n\right)^{-2}=\sum_{n=1}^{\infty}\left(2n-1\right)^{-2}+\frac{1}{4}\sum_{n=1}^{\infty}n^{-2}.
\]
\end{rem}
\begin{prop}
\label{prop:brownian}
In the above hypotheses, given $r$ in the range of $g$ such that
\begin{align}
\label{cond_r_alpha_beta}
\alpha_2r\leq -\left(\beta_{2}+\frac{\| g''\|_{C_b(\R)}}{2}\right), 
\end{align}
we define 
\begin{align*}
G(x)=\int_{0}^{1}g(x(s))\ ds-r, \quad x \in X.
\end{align*}
Then, $G$ satisfies Hypotheses \ref{claim:regularity} and \ref{claim:Gaussian-curvature}.
\end{prop}
\begin{proof}
To prove that Hypothesis \ref{claim:regularity} is satisfied it is enough to argue as in \cite[Proposition 5.1]{add2}). 

We stress that $G$ is $H$-differentiable at every $x\in X$ and, for every $h\in H$,
\begin{align}
\label{def_der_G_ex}
\langle \nabla_{H}G(x),h\rangle _{H}=\int_{0}^{1}g'(x(s))h(s)\ ds.
\end{align}
This implies that, for every $x\in X$, from \eqref{eq:asymptote-1} and the definition of $h_1$ it follows that 
\begin{align*}
\|\pi_1(\nabla_HG(x))\|_H^2
=  & \left(\int_0^1g'(x(s))h_1(s)ds\right)^2\geq \frac{8a}{\pi^2}.
\end{align*}
Hence, the first part of Hypothesis \ref{claim:Gaussian-curvature} is satisfied with $n_0=1$.

Moreover, $D_{H}G$ is $H$-differentiable at every $x\in X$ and, for every $h,k\in H$,
\[
\langle D_{H}^{2}G(x)(h),k\rangle _{H}=\int_{0}^{1}g''(x(s))h(s)k(s)\ ds,
\]
and so, fixed $h,k\in H$, for every $x_{1},x_{2}\in X$ we have
\begin{align*}
|\langle D_{H}^{2}G(x_{1})(h),k\rangle _{H}-\langle D_{H}^{2}G(x_{2})(h),k\rangle _{H}|
\leq & |\int_{0}^{1}(g''(x_{1}(s))-g''(x_{2}(s)))h(s)k(s)\ ds| \\
\leq & [g'']_\alpha\int_{0}^{1}|x_{1}(s)-x_{2}(s)|^{\alpha}|h(s)||k(s)|\ ds \\
\leq &  [g'']_\alpha\|h\|_{H}\|k\|_{H}\|x_{1}-x_{2}\|_{X}^{\alpha},
\end{align*}
where $[g'']_\alpha$ is the H\"older seminorm of $g''$. It follows that $G$ is a $C^{2,\alpha}$ function on every subspace $F$ of $H$ with $\mbox{dim}(F)<\infty$, which means that $\mathcal O_n$ is a $C^{2,\alpha}$-regular open set in $\R^n$ for every $n\in\N$. 

It remains to show that $H_{\partial O_n}^\gamma\geq 0$ on $\partial O_n$ for every $n\in\N$. To this aim, for every $s\in[0,1]$ we set
\begin{align*}
f(s)
= & \sum_{n=1}^{\infty} (h_n(s))^2
= \sum_{n=1}^{\infty}2\pi^{-2}\left(n-\frac12\right)^{-2}\left(\sin\left(\Big(n-\frac12)\Big)\pi s\right)\right)^2 \\
= & \frac12-\sum_{n=1}^{\infty}\pi^{-2}\left(n-\frac12\right)^{-2}\cos\big((2n-1)\pi s\big),
\end{align*}
where the last equality follows from \eqref{eq:basel_modified}. Since 
\begin{align*}
\int_{0}^{1}\left(\cos(2n-1)\pi s\right)ds=0, \qquad n\in\N
\end{align*}
and the series which defines $f$ totally converges, we conclude that
\begin{align*}
\int_{0}^{1}f(s)\ ds=\frac{1}{2}.
\end{align*}

Let us fix $m\in\N$ and $x\in\partial O_{m}$, and let us set $\varphi_{m}=\nabla_{H}G_{m}(x)=\pi_m(\nabla_HG(\pi_m(x)))$. It follows that $\varphi_{m}\in H$ and $\|\varphi_{m}\|_{H}\leq\| g\| _{C_{b}^{1}(\R)}$. We introduce the function
\begin{align*}
f_{m}(s) & =\sum_{n=1}^{m}(h_n(s))^2,
\end{align*}
Then, $f_{m}>0$ on $[0,1]$, $\|f_m\|_{L^1(0,1)}\leq \|f\|_{L^1(0,1)}= \frac12$ and
\begin{align}
\label{ex_traccia_D2G}
{\rm Tr}_H[D^2_HG_m(x)]
= \sum_{n=1}^m\int_0^1g''((\pi_mx)(s))(h_n(s))^2ds =\int_0^1g''((\pi_mx)(s))f_m(s)ds .
\end{align}
From \eqref{def_der_G_ex} and \eqref{ex_traccia_D2G} it follows that
\begin{align*}
LG_{m}(x)=\int_{0}^{1}g''((\pi_{m}x)(s))f_{m}(s)ds-\int_{0}^{1}g'((\pi_{m}x)(s))(\pi_{m}x)(s)ds,
\end{align*}
and so
\begin{align}
& H_{\partial O_m}^\gamma(x)= \notag \\
= & \frac{LG_{m}(x)}{\|\nabla_{H}G_{m}(x)\|_{H}}-\frac{\langle D_{H}^{2}G_{m}(x)\nabla_{H}G_{m}(x),\nabla_{H}G_{m}(x)\rangle_H}{\|\nabla_{H}G_{m}(x)\|_{H}^{3}} \notag \\
= & \frac{\int_{0}^{1}g''((\pi_{m}x)(s))f_{m}(s)ds-\int_{0}^{1}g'((\pi_{m}x)(s))(\pi_{m}x)(s)ds}{\|\varphi_{m}\|_{H}} 
-\frac{\int_{0}^{1}g''((\pi_{m}x)(s))(\varphi_{m}(s))^2ds}{\|\varphi_{m}\|_{H}^{3}} \notag \\
= & \|\varphi_{m}\|_{H}^{-1}\Big(\int_{0}^{1}g''((\pi_{m}x)(s))\left(f_{m}(s)-\frac{(\varphi_{m}(s))^2}{\|\varphi_{m}\|_{H}^{2}}\right)ds 
-\int_{0}^{1}g'((\pi_{m}x)(s))(\pi_{m}x)(s)ds\Big).
\label{ex_2_Hmgamma}
\end{align}
An explicit computation on $\varphi_m$ gives
\begin{align}
\label{ex_2_stima_rapp_phi}
(\varphi_m(s))^2
= & \left(\sum_{n=1}^m\langle \nabla_HG_m(x),h_n\rangle_H h_n(s)\right)^2
\leq 
\|\varphi_m\|_H^2f_m(s), \quad s\in[0,1].
\end{align}
Further, from \eqref{eq:asymptote-1} we infer that
\begin{align}
\label{ex_2_stima_der_Hmgamma}
\int_{0}^{1}g'((\pi_{m}x)(s))(\pi_{m}x)(s)\ ds
\leq & \alpha_2\int_0^1 g((\pi_mx)(s))ds+\beta_2,
\end{align}
and so from \eqref{ex_2_Hmgamma}, \eqref{ex_2_stima_rapp_phi} and \eqref{ex_2_stima_der_Hmgamma} we get
\begin{align*}
 & H_{\partial O_m}^\gamma(x) \geq  \\
& \geq  \|\varphi_m\|_H^{-1}
\Bigg(\int_{\{s\in(0,1):g''((\pi_mx)(s))\leq 0\}}g''((\pi_mx)(s))\left(f_m(s)-\frac{(\varphi_m(s))^2}{\|\varphi_m\|_H^2}\right)ds  \\
&\quad - \alpha_2\int_0^1 g((\pi_mx)(s))ds-\beta_2\Bigg) \\
& \geq  \|\varphi_m\|_H^{-1}
\Bigg(\int_{\{s\in(0,1):g''((\pi_mx)(s))\leq 0\}}g''((\pi_mx)(s))f_m(s)ds -\alpha_2 \int_0^1 g((\pi_mx)(s))ds-\beta_2\Bigg) \\
& \geq  \|\varphi_m\|_H^{-1}
\Bigg(-\|g''\|_{C_b(\R)}\int_0^1f_m(s)ds -\alpha_2(G_m(x)+r)-\beta_2\Bigg) \\
& =   \|\varphi_m\|_H^{-1}
\Bigg(-\frac{\|g''\|_{L^\infty(\R)}}2-\alpha_2 r-\beta_2\Bigg),
\end{align*}
where we have used the fact that $G_m(x)=0$ for $x\in\partial O_m$. Finally, from \eqref{cond_r_alpha_beta} it follows that
\begin{align*}
-\frac{\|g''\|_{L^\infty(\R)}}2-\alpha_2 r-\beta_2\geq 0,
\end{align*}
which gives $\mathcal H_m^\gamma(x)\geq0$ for every $x\in \partial O_m$, for every $m\in \N$. This implies that also the last point in Hypothesis \ref{claim:Gaussian-curvature} is fulfilled.
\end{proof}

\subsubsection{Brownian bridge on $0$}

We consider a pinned Wiener space, which describes a Brownian bridge with starting point at $0$ and subjected to the condition that at $1$ the arrival point is $0$. In this setting, $X=L^{2}[0,1]$, the Cameron-Martin
space is $H=H_{0}^{1}(0,1)$ and $\gamma$ is the pinned Wiener measure on $X$, see \cite{Revuz}. We recall that it is defined an orthonormal basis $\{e_{n}=\sqrt 2\sin(\pi n\cdot):n\in\N\}$ of $X$ with eigenvalues $\lambda_{n}=(\pi n)^{-2}$ with respect to the covariance operator $Q$ of $\gamma$, and $\{\sqrt{\lambda_n}e_n=\sqrt{2}\pi^{-1}n^{-1}\sin(n\pi\cdot):n\in\N\}$ is an orthonormal basis of $H$ of eigenvectors of $Q$.

We introduce a function $g$ which satisfies the same assumptions as in Example \ref{ex:br_mot}. Arguing as in the proof of Proposition \ref{prop:brownian} we can prove the following result. 
\begin{prop}
In the above hypotheses, given $r$ in the range of $g$ and 
\begin{align*}
\alpha_2r\leq -\left(\beta_{2}+\frac{\| g''\|_{C_b(\R)}}{6}\right), 
\end{align*}
we define
\[
G(x)=\int_{0}^{1}g(x(s))\ ds-r, \quad x\in X.
\]
Then, $O=G^{-1}((-\infty,0))$ satisfies Hypotheses \ref{claim:regularity} and \ref{claim:Gaussian-curvature}.
\end{prop}
The main difference in the proof is that, in this case, we get $f=s-s^{2}$ and $f_{m}=2\sum_{n=1}^{m}\lambda_{n}\sin(n\pi s)^{2}$ (see \cite[Example 5.4]{Dap} for details).

\vspace{3mm}
{\bf Acknowledgments.} The authors are grateful to Michael R\"ockner for the suggestion to extend to the finite dimensional setting the result in \cite[Appendix 8]{Bar}.

\end{document}